\documentclass[pdflatex,sn-mathphys-num]{sn-jnl}

\usepackage{graphicx}%
\usepackage{multirow}%
\usepackage{amsmath,amssymb,amsfonts}%
\usepackage{amsthm}%
\usepackage{mathrsfs}%
\usepackage[title]{appendix}%
\usepackage[table]{xcolor}%
\usepackage{textcomp}%
\usepackage{manyfoot}%
\usepackage{booktabs}%
\usepackage{algorithm}%
\usepackage{algorithmicx}%
\usepackage{algpseudocode}%
\usepackage{listings}%
\usepackage{todonotes}
\usepackage{siunitx}
\usepackage{a4wide}

\theoremstyle{thmstyleone}%
\theoremstyle{thmstyletwo}%

\theoremstyle{thmstylethree}%

\newcommand{\skw}[1]{\hat{#1}}

\newcommand{\norm}[1]{\left\Vert #1 \right\Vert}

\newcommand{\reals}{\mathbb R}
\newcommand{\liegroup}{SO(3)}
\newcommand{\liealgebra}{\mathfrak{so}(3)}

\newcommand{\domain}{\Omega_l}

\newcommand{\av}{\boldsymbol a}
\newcommand{\bv}{\boldsymbol b}

\newcommand{\ev}{\boldsymbol e}
\newcommand{\ex}{\ev_1}
\newcommand{\ey}{\ev_2}
\newcommand{\ez}{\ev_3}

\newcommand{\gv}{\boldsymbol g}
\newcommand{\uv}{\boldsymbol u}

\newcommand{\I}{\mathbf I}
\newcommand{\zerov}{\boldsymbol 0}

\newcommand{\rv}{\boldsymbol r}
\newcommand{\psiv}{\boldsymbol \psi}

\newcommand{\rotmat}{\boldsymbol \Lambda}
\newcommand{\rot}{\rotmat}
\newcommand{\rotvec}{\psiv}
\newcommand{\rotangle}{\psi}
\newcommand{\tangent}{\mathbf T}

\newcommand{\kappav}{\boldsymbol{\kappa}}
\newcommand{\gammav}{\boldsymbol{\gamma}}
\newcommand{\thetav}{\boldsymbol{\theta}}
\newcommand{\mv}{\boldsymbol{m}}
\newcommand{\nv}{\boldsymbol{n}}
\newcommand{\Mv}{\boldsymbol{M}}
\newcommand{\Nv}{\boldsymbol{N}}
\newcommand{\qv}{\boldsymbol{q}}
\newcommand{\cv}{\boldsymbol{c}}

\newcommand{\Wint}{W_{int}}
\newcommand{\Wext}{W_{ext}}

\newcommand{\Cgamma}{\mathbb C_\gamma}
\newcommand{\Ckappa}{\mathbb C_\kappa}

\newcommand{\EA}{EA}
\newcommand{\GAy}{GA_2}
\newcommand{\GAz}{GA_3}
\newcommand{\GIt}{GI_t}
\newcommand{\EIy}{EI_2}
\newcommand{\EIz}{EI_3}

\newcommand{\err}{e}
\newcommand{\nel}{n_{\rm elem}}

\newcommand{\Fv}{\boldsymbol F}

\graphicspath{{./figures/}}

\usepackage{etoolbox}
\makeatletter
\patchcmd{\@maketitle}
  {{\printabstract\par}}%
  {{\ifx\@motto\empty\else\@motto\fi\printabstract\par}}%
  {}{\typeout{sn-jnl.cls: Failed to patch \string\@maketitle}}
\makeatother

\begin{document}

\title[Mixed Finite Elements for Geometrically Exact Beams]{Mixed Finite Elements for Geometrically Exact Beams using Discontinuous Rotations and Discrete Curvature}


\author*[1]{\fnm{Alexander} \sur{Humer}}\email{alexander.humer@jku.at}

\author[2]{\fnm{Ivo} \sur{Steinbrecher}}\email{ivo.steinbrecher@unibw.de}

\author[1]{\fnm{Astrid} \sur{Pechstein}}\email{astrid.pechstein@jku.at}

\affil*[1]{\orgdiv{Institute of Technical Mechanics}, \orgname{Johannes Kepler University Linz}, \\ \orgaddress{\street{Altenberger Str. 69}, \postcode{4040} \city{Linz}, \country{Austria}}}

\affil[2]{\orgdiv{Institute for Mathematics and Computer-Based Simulation}, \orgname{University of the Bundeswehr Munich}, \orgaddress{\street{Werner-Heisenberg-Weg 39}, \postcode{85579} \city{Neubiberg}, \country{Germany}}}


\motto{\centering Dedicated to Professor Hans Irschik on the occasion of his 75th birthday}

\abstract{
    We propose a novel mixed finite-element formulation for geometrically exact (Simo--Reissner) beams that introduces the moment vector as additional independent field. 
    The specific mixed form allows for an element-local, discontinuous approximation of rotations, which is key to a simple and efficient discretization framework.
    The concept of discrete curvature provides a mathematically consistent treatment of rotation discontinuities.
    For linear constitutive laws, the mixed form is derived via a Legendre transform of the curvature-related strain energy.
    Objectivity is retained at the discrete level by interpolating relative rotations through a multiplicative split of the rotation field; path-independence is inherent to the total Lagrangian setting and verified numerically.
    Several benchmarks demonstrate optimal rates of convergence and accuracy, irrespective of the beam's slenderness and order of approximation.
    Notably, the lowest-order element entirely avoids rotation interpolation by employing element-constant rotations only.
}

\maketitle

\section{Introduction}
\label{sec:intro}
The motivation for our research is relatively simple: finite rotations are an intriguing and, at times, frustrating subject, which, in combination with structural theories, are the perfect playground for those interested in numerical methods, in particular. 
And this still holds to this day, even though it is certainly no exaggeration to say that beam theories---and associated numerical methods---have come a long way. 
It has been nearly 300 years since Euler, in the appendix to his book~\cite{Euler1744}, derived solutions for the elastic curve (``curvae elasticae'') and established what is now known as Euler's elastica.
Subsequent developments are associated with such great names as Kirchhoff~\cite{Kirchhoff1859} and Timoshenko~\cite{Timoshenko1921}, whose contributions---among several others---were key for what we now call the \emph{geometrically exact beam theory}, which translates Timoshenko's kinematic assumptions to problems involving large deformations.
The geometrically exact beam theory is commonly attributed to the likes of Reissner~\cite{Reissner1972,Reissner1973,Reissner1981}, Antman~\cite{Antman1973} and Simo~\cite{Simo1985,Simo1986,Simo1986a,Simo1986b}.
Due to the ambiguity of what counts as ``geometrically exact'', the term ``Simo--Reissner beam'' theory is frequently used to differentiate from large deformation ``Kirchhoff--Love beams'', which do not introduce kinematic approximations either, but are rigid in shear.
The term also acknowledges the theoretical foundation established by Reissner in terms of generalized strain measures, which are derived from variational considerations---first for the planar case~\cite{Reissner1972}, later for general spatial deformation~\cite{Reissner1973, Reissner1981}.
For the planar case, Irschik and Gerstmayr~\cite{Irschik2009,Irschik2011} established a continuum-mechanics interpretation of Reissner's strain measures and stress-resultants, based on which an extension to thermoelasticity was presented in~\cite{Humer2021}.

Though closed-form solutions in terms of elliptic integrals/functions exist for selected static problems with concentrated loads~\cite{Humer2013, Humer2019}, we generally need to resort to numerical methods to construct approximate solutions.
Simo and Vu-Quoc were arguably the first to present a successful finite-element formulation for geometrically exact beams, this time in reverse order, i.e., the spatial case came first~\cite{Simo1985,Simo1986}, then the planar formulation~\cite{Simo1986a,Simo1986b}.
Their contributions, which continue to influence the numerical analysis of beams to this day, accounted for the Lie-group structure of the configuration space emanating from the necessity to describe large rotations of cross-sections in geometrically exact beams. 
We do not want to conceal that there were successful formulations for nonlinear beam problems before; we mention Bathe and Bolourchi~\cite{Bathe1979} as an example, since one of their benchmarks is also considered herein. 
In the wake of Simo’s contributions, however, developments in the field of finite beam elements really gathered momentum.
Subsequent works focused primarily on how to best handle large rotations. 
In addition to the parameterization of rotations themselves, both formulations based on absolute rotations, see, e.g.,~\cite{Cardona1988, Ibrahimbegovic1995}, and incremental formulations~\cite{ibrahimbegovic1997a}, for which the question of path-independence is particularly relevant, were introduced.
The review paper by Meier et al.~\cite{Meier2019} gives an excellent overview of the developments in this field.

The contributions of Crisfield and Jelenić~\cite{Crisfield1999, Jelenic1999} marked an important turning point, which is also relevant for the present work.
They demonstrated that discretization schemes that had been used up to this point failed to preserve the objectivity of generalized strain measures of the geometrically exact beam theory.
As a remedy, they proposed to interpolate relative local rotations rather than total rotation vectors.
Unit quaternions, which are also referred to as Euler parameters~\cite{Wehage1984,Nikravesh1985}, are an alternative way to represent rotations by means of four parameters, which was adopted in the context of geometrically exact beams, e.g., by Zupan et al.~\cite{Zupan2009, Zupan2013}, and, more recently, by Wasmer and Betsch~\cite{Wasmer2024}.
In the present work, we follow Crisfield and Jelenić~\cite{Crisfield1999} by using rotation vectors for the parameterization of rotations.
Special care is taken to ensure the representation of rotations through shape functions of arbitrary order in an objective manner by introducing a multiplicative split in a similar manner.
In these implementations, we make use of quaternion algebra on the algorithmic level to efficiently evaluate composite rotations.
Over the last decade, formulations based on the concept of the \emph{special Euclidean group} $SE(3)$ have gained popularity~\cite{Sonneville2014,Sonneville2017,harsch2023}. 
This so-called \emph{motion formalism} is based on the concept of jointly representing both translations and rotations of cross-sections by a single Lie group, whose elements are referred to as ``motions'', rather than treating displacements and rotations separately. 
In~\cite{Sonneville2014,harsch2023} a two-noded element was presented, where the interpolation of  motions was realized by a linear scaling of the relative twist between the nodal values, in analogy to what Crisfield and Jelenić had proposed for rotations~\cite{Crisfield1999}.
To avoid any mesh-related bias due to the selection of a reference node, Sonneville et al.~\cite{Sonneville2017} proposed an implicit scheme, which allows for Lagrange polynomials of arbitrary order to be used for the interpolation of relative motions. 
As implicit interpolation departs from conventions underlying many finite-element frameworks, the algorithmic implications of such approaches are considerable.
Due to the coupled representation of translations and rotations, beam elements based on the motion formalism have been shown to be less susceptible to \emph{locking phenomena}, i.e., excessive stiffness originating from a spurious coupling among strain measures.
Reduced integration is known as an effective measure to alleviate locking in various kinds of problems, which, in the context of beams, was already applied in Simo and Vu-Quoc's original beam element~\cite{Simo1986}.
\emph{Mixed methods} offer another way to address locking.
Although mixed methods are not as common for beams as they are in shell formulations, the literature does contain several approaches employing two- and three-field formulations based on variants of the Hellinger--Reissner~\cite{Noor1981, Nukala2004, Herrmann2026} and Hu--Washizu~\cite{Saje1990,Saje1991} variational principles, where cross-sectional stress resultants, in case of the latter, generalized strains complement displacements and rotations as independent fields.
In a recent contribution, Herrmann et al.~\cite{Herrmann2026} proposed a Hellinger--Reissner-type formulation for geometrically exact beams using quaternions for the parameterization of rotations.
They showed that problems of Kirchhoff--Love beams are naturally contained as a special case of a vanishing shear compliance.

The mixed formulation for geometrically exact beams presented in this work departs from previous approaches in many respects. 
Our goal is to derive an equally simple and effective mixed approach that meets the essential requirements of objectivity and path-independence.
For this purpose, we introduce only the moment vector as an additional independent field variable, while retaining the primal form with respect to forces and displacements.
We draw some of our inspiration from mixed methods for Kirchhoff--Love shells~\cite{neunteufel2019,neunteufel2021,Platzer2025a,Platzer2025b}, which share our notion of \emph{discrete curvature}.
%
%
The discrete curvature allows for a discontinuous interpolation of rotations, which (largely) eliminates locking phenomena.
Being a total Lagrangian approach that forgoes any increments in rotational variables, the proposed formulation is a priori path-independent. 
Objectivity of discrete strains is established by interpolating relative rotations, which are recovered ``for free'' upon a multiplicative split of rotations into constant and higher-order contributions that is realized at shape-function level.
By including nodal rotations as additional variables, we are able to impose boundary conditions just as in conventional irreducible approaches.
These ingredients allow for a locking-free \emph{lowest-order element} that uses element-wise constant rotations, thereby avoiding any interpolation of rotations.
Our formulation naturally extends to arbitrary higher orders, for which we can recover ideal convergence rates.
Moreover, we show that the formulation can be immediately applied to structures with kinks and branches, which are common features in real-world problems. 


The present paper is organized as follows: In Section~\ref{sec:beam_theory}, we briefly revisit the geometrically exact beam theory, which provides the theoretical foundation of the proposed finite-element formulation.
We start from the kinematics of beams being represented as material Cosserat lines, in which rotations and their parameterization play a central role. 
Strain measures of the theory are discussed, and the equilibrium relations are provided in both strong and weak form.
Based on the latter, Section~\ref{sec:fe} establishes the proposed mixed formulation, in which the moment vector serves as independent field.
The notion of discrete curvature, which is key to the discontinuous interpolation of rotations, is introduced subsequently.
In Section~\ref{sec:examples}, we apply the proposed beam element in several numerical examples to demonstrate its effectiveness by means of convergence rates. 
Additionally, we numerically verify that both objectivity and path-independence of the formulation are preserved upon discretization.
Section~\ref{sec:conlusion} concludes the paper with a brief discussion of the main contributions and an outlook on possible extensions.

\section{Geometrically exact beam theory}
\label{sec:beam_theory}

\subsection{Kinematics of shear-deformable beams}
\begin{figure}[htb!]
    \centering
    \includegraphics[scale=0.8]{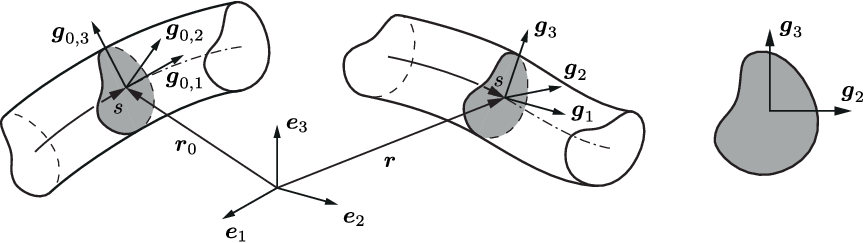}
    \caption{Kinematics of a geometrically exact beam: The beam's undeformed (indicated by subscripts ``$0$'') and deformed configurations are determined by the centerline position $\rv$ and the orientation of the cross-section, which is represented by the triad $\{\gv_1,\gv_2,\gv_3\}$ (adapted from~\cite{Steinbrecher2020}).}
    \label{fig:kinematics}
\end{figure}

%
In what follows, we briefly review the essential constituents of what is usually referred to as geometrically exact beam theory.
We assume cross-sections of a beam to remain plane and undistorted in the course of deformation, but, in presence of shear deformation, cross-sections do not necessarily stay perpendicular to the beam's centroid axis (or centerline).
As a consequence, the beam's configuration is fully determined by the position of its centerline and the orientation of its cross-sections, see Fig.~\ref{fig:kinematics} for an illustration.
From a modeling point of view, beams can be regarded as material \emph{Cosserat lines}, i.e., material points have both translational and rotational degrees of freedom.

More specifically, let $\rv_0$ denote the position of the beam's axis in the undeformed configuration, in which it describes a spatial curve that is parameterized by the arclength coordinate $s \in \domain$: 
\begin{align}
    &s \mapsto \rv_0 \in \reals^3 . 
\end{align}
In the undeformed configuration, which is indicated by subscripts `0', the orientation of cross-sections relative to some global Cartesian frame $\left\{ \ex, \ey, \ez \right\}$ is represented by the triad $\rot_0$, 
\begin{align}
    s \mapsto \rot_0 = \gv_{0,i} \otimes \ev_i \in \liegroup ,
\end{align}
i.e., rotation tensors belonging to the group $\liegroup$ of orthogonal transformations.
By convention, $\gv_{0,1}$ is tangential to the beam's (undeformed) axis, which implies $\gv_{0,1} = \rv_0^\prime := \partial \rv_0 / \partial s$; unit vectors $\gv_{0,2}$ and $\gv_{0,3}$ span the cross-section in the undeformed state.

Under the action of external loads, the beam's axis is displaced and cross-sections rotate. 
As for the undeformed state, the deformed (current) configuration of the beam is described by the axis' position and the orientation of cross-sections (with respect to the global frame), which we subsequently denote by
\begin{align}
    &s \mapsto \left( \rv, \rot \right) \in \reals^3 \times \liegroup,
    &\rot = \gv_i \otimes \ev_i .
\end{align}
Note that, owing to shear deformation, the unit vector $\gv_1$ generally is no longer tangential to the beam's axis in the deformed configuration.
The remaining elements of the triad $\gv_2$ and $\gv_3$ still span the cross-section.

\subsection{Rotations and rotation vectors}
Undoubtedly, rotations lie at the heart of most, if not all, major challenges that come along with numerical methods for Cosserat-like continua involving the orientation of material points, lines or cross-sections as an independent field.
Due to the specific nature of the special orthogonal group $\liegroup$, the configuration space of geometrically exact beams constitutes a nonlinear manifold that is endowed with the properties of Lie groups.
The non-additivity of elements of the Lie group and the non-commutativity of the group product contrast 
our intuition of interpolation in vector spaces.

A first question to be answered is which parameterization of rotations is to be chosen.
Among the nine components of an orthogonal tensor $\rot$, only three are independent due to the  defining property of $\liegroup$, i.e.,
\begin{equation}
    \liegroup = \left\{ \rot \in \reals^{3 \times 3} \; | \; \rot^{-1} = \rot^T, \; \det \rot = 1 \right\} .
\end{equation}
We refer to \cite{shabana1998b,Geradin2001,Crisfield2003} for a discussion of the various possibilities to parameterize rotation tensors as, e.g., Euler angles, Rodrigues parameters or unit quaternions, to mention some. 
To this end, we use the concept of rotation vectors, which have an intuitive geometric interpretation. The orientation of a rotation vector $\rotvec \in \reals^3$, i.e., $\ev_\psi = \rotvec / \norm{\rotvec}$ coincides with the axis of rotation; the scalar-valued angle of rotation about this axis is given by its length $\psi = \norm{\rotvec}$.
The \emph{exponential map} provides a smooth mapping from the linear space of rotation vectors to nonlinear manifold of orthogonal tensors,
\begin{align} \label{eq:exp_so3}
    \rot &= \exp \skw \rotvec = \sum_{k=0}^\infty \frac{1}{k!} \bigl( \skw\rotvec \bigr)^k ,
\end{align}
where a hat symbol has been introduced to indicate skew-symmetric tensors isomorphic to three-dimensional vectors in the sense that 
\begin{align} 
    &\skw \av \cdot \bv = \av \times \bv ,
    &\skw \av = - \skw \av^T ,  &
    &\av, \bv \in \reals^3 .
\end{align}
For the matrix exponential of skew-symmetric tensors, \emph{Rodrigues' formula} provides us with a closed-form relation to compute the rotation tensor for a given rotation vector,
\begin{equation}
    \rot = \I + \frac{\sin \rotangle}{\rotangle} \skw \rotvec + \frac{1 - \cos \rotangle}{\rotangle^2} (\skw \rotvec)^2 .
\end{equation}
From an algebraic point of view, the space of skew-symmetric tensors $\skw \rotvec$ constitutes the \emph{Lie algebra} $\liealgebra$, i.e., the tangent space of the Lie group $\liegroup$ at the identity. 
The exponential map, $\exp : \liealgebra \to \liegroup$, is a surjective map from the Lie algebra to the Lie group.
It is, however, not injective, due to the fact that parallel rotation vectors with lengths offset by multiples of $2\pi$ represent the same rotation.
We nevertheless define its right inverse yielding a rotation vector $\sphericalangle({\rot})$ for a given rotation $\rot$, which is bounded in absolute value by $[0,\pi]$ and satisfies
\begin{align}
    \rot = \exp(\skw{\psiv}) \qquad \text{for } \psiv = \sphericalangle({\rot}).
\end{align}

Successive rotations, e.g., $\rot_1 \in \liegroup$ followed by $\rot_2 \in \liegroup$, are represented by the group product, i.e., matrix multiplication,
\begin{equation}
    \rot = \rot_2 \cdot \rot_1 = \exp \skw \rotvec_2 \cdot \exp \skw \rotvec_1 \in \liegroup .
\end{equation}
%

Care must be taken when employing \emph{variations} of rotation vectors or tensors.
One can express the variation of the rotation in terms of the an additive perturbation of the rotation vector $\delta \rotvec$ as the directional derivative of the exponential map, which, in view of the chosen parameterization of rotations, naturally plays an important role within the proposed approach:
\begin{align}
    &\delta \rot = \lim_{\epsilon \to 0}\frac{d}{d \epsilon} \exp \bigl( \skw\rotvec + \epsilon \delta \skw \rotvec \bigr) = \frac{\partial \rot}{\partial \rotvec}  \cdot \delta \rotvec .
\end{align}
Alternatively, one can regard the variation $\delta \rot$ as a (linearized) consecutive rotation following $\rot$. The according \emph{spin vector} is commonly referred to as $\delta \thetav$, and is obtained by means of the directional derivative as
%
\begin{align}
    &\delta \rot = \lim_{\epsilon \to 0}\frac{d}{d \epsilon} \exp \bigl( \epsilon \delta \skw \thetav \bigr) \cdot \rot = \delta \skw \thetav \cdot \rot 
    & \leftrightarrow &
    & \delta \skw \thetav = \delta \rot \cdot \rot^T \in \liealgebra .
\end{align}
From the above representations of $\delta \rot$, the relation between the spin vector and the variation of the rotation vector can be established in terms of the tangent map (or application) $\tangent (\rotvec)$ as
\begin{equation}
    \delta \thetav = \tangent \cdot \delta \rotvec , 
\end{equation}
for which we have the closed-form expression
\begin{equation} \label{eq:tangent_so3}
    \tangent = \I + \frac{1 - \cos \rotangle}{\rotangle^2} \skw \rotvec + \frac{1}{\rotangle^2} \left( 1 - \frac{\sin \rotangle}{\rotangle} \right) (\skw \rotvec)^2 .
\end{equation}

    Although both the exponential map~\eqref{eq:exp_so3} and the tangent map~\eqref{eq:tangent_so3} are well-defined for the limiting case $\rotangle \rightarrow 0$, we use the respective second-order series expansions in our implementation for small angles $\rotangle < \epsilon_\rotangle$ in order to avoid numerical instabilities: 
    \begin{align}
        &\rotmat \approx \I + \skw \rotvec + \frac 1 2 (\skw \rotvec)^2 , &
        &\tangent \approx \I + \frac 1 2 \skw \rotvec + \frac 1 6 (\skw \rotvec)^2 .
    \end{align}

\subsection{Strain measures and objectivity}
To quantify the state of deformation, we need appropriate strain measures that obey the fundamental requirement of \emph{objectivity}, i.e., invariance under arbitrary rigid-body motions that are superimposed to the deformed configuration. 
In his fundamental contributions, Reissner~\cite{Reissner1972,Reissner1973,Reissner1981} motivated the definitions of what he referred to as force and moment strains based on the equivalence of local balance relations of linear and angular momentum and a postulated structure of the principle of virtual work.
The strain measures are meant to quantify (local) changes in the configuration in the course of deformation.
In what follows, we employ the material representation of strains, in which the difference between the current configuration and the undeformed state is expressed in the common global frame.
In this sense, the (material) force strain $\gammav$ related to shear and extension is defined as
\begin{equation}
    \label{eq:gammav}
    \gammav = \rot^T \cdot \rv^\prime - \rot_0^T \cdot \rv_0^\prime = \rot^T \cdot \rv^\prime - \ex .
\end{equation}
A change in twist and curvature of the beam is quantified by the (material) moment strain $\kappav$, which is given by
\begin{equation}
    \skw \kappav = \rot^T \cdot \rot^\prime - \rot_0^T \cdot \rot_0^\prime .
    \label{eq:define_kappa}
\end{equation}
The moment strain, which we also refer to as curvature subsequently, can also be expressed in terms of the tangent map~\eqref{eq:tangent_so3} as 
\begin{align}
    \kappav 
    &= \tangent^T(\rotmat) \cdot \rotvec^\prime - \tangent^T(\rotmat_0) \cdot \rotvec_0^\prime . 
\end{align}
%
    The above definitions of strains corresponds to the idea of a fictitious reference configuration, in which the beam's axis is straight and aligned with the unit normal $\ex$ of the global frame.
    In case of structures with kinks (``slope discontinuities'') and branches, however, such reference configuration may be somewhat difficult to visualize. 
    The proposed formulation, though, naturally extends to problems characterized by geometric features of this kind.

%
Let $(\rv_{rb}, \rot_{rb}) \in \reals^3 \times \liegroup$ denote a \emph{constant} translation and rotation, which we superimpose the beam's current configuration. 
Upon such rigid-body motion, we obtain the configuration $(\bar \rv, \bar \rot)$, 
\begin{align}
    \bar \rv &= \rot_{rb} \cdot \left( \rv_{rb} + \rv \right) , &
    \bar \rot &= \rot_{rb} \cdot \rot ,
\end{align}
for which we can immediately convince ourselves by direct computation that the strain measures remain unaffected,
\begin{align}
    \bar \gammav &= \rot^T \cdot \rot_{rb}^T \cdot \rot_{rb} \cdot \rv^\prime - \ex = \gammav , &
    \skw{\bar{\kappav}} &= \rot^T \cdot \rot_{rb}^T \cdot \rot_{rb} \cdot \rot^\prime - \rot_0^T \cdot \rot_0^\prime = \skw\kappav ,
\end{align}
i.e., they are \emph{objective}.
In their seminal paper, Crisfield and Jelenić~\cite{Crisfield1999,Jelenic1999} have shown that property of objectivity generally is not preserved upon a discretization of the rotation field.
In particular, the pioneering formulations of Simo and Vu-Quoc~\cite{Simo1985,Simo1986} and subsequent contributions based upon these fail to meet this requirement, which is typically considered essential in continuum theories and discrete approximations thereof.
The finite-element formulation proposed subsequently provides an elegant way to retain the objectivity of strains upon a spatial discretization.

\subsection{Equilibrium relations and weak form}
Prior to the finite-element formulation, we recall the relations that govern the deformation of the beam under external loads.
Staying with the static case, the local equilibrium relations of a Cosserat line are given by
\begin{align}
    &\nv^\prime + \qv = \zerov , 
    &\mv^\prime + \rv^\prime \times \nv + \cv = \zerov ,
\end{align}
where $\nv$ and $\mv$ represent the cross-sectional forces and moments in the \emph{current configuration}; externally imposed distributed forces and couples are denoted by $\qv$ and $\cv$, respectively.
A pull-back with the rotation $\rot$ gives the corresponding material forces and moments
\begin{align} \label{eq:material_stress_resultants}
    &\Nv = \rot^T \cdot \nv, & 
    &\Mv = \rot^T \cdot \mv .
\end{align}
%
The variational (or weak) formulation is derived considering admissible virtual displacements $\delta \rv$ and rotations $\delta \thetav$.
From a mechanical perspective, the weak form constitutes what is referred to as the \emph{principle of virtual work}, which states that the virtual work performed by internal forces upon a virtual deformation characterized through $\delta \rv$ and $\delta \thetav$ from the current configuration equals the virtual work performed by the external loads, i.e., $\delta \Wint = \delta \Wext$, which are given by 
\begin{equation}
    \delta \Wint = \int_{\domain} \left( \Nv \cdot \delta \gammav + \Mv \cdot \delta \kappav \right) ds ,
    \label{eq:delta_W_int}
\end{equation} 
and
\begin{equation}
    \delta \Wext = \int_{\domain} \left( \qv \cdot \delta \rv + \cv \cdot \delta \thetav \right) ds
    + \left[ \nv \cdot \delta \rv + \mv \cdot \delta \thetav \right]_0^L ,
\end{equation}
respectively.

For elastic materials, the (material) stress resultants $\Nv$ and $\Mv$ depend on the current strain measures $\gammav$ and $\kappav$ only. 
In the following, we assume the existence of strain energies densities $\varphi_\gamma$ and $\varphi_\kappa$ describing these relations through
\begin{align}
    \label{eq:constitutive_equations}
    &\Nv = \frac{\partial \varphi_\gamma}{\partial \gammav}, &
    &\Mv = \frac{\partial \varphi_\kappa}{\partial \kappav}.
\end{align}
In the case of linear constitutive relations, the second-order tensors $\Cgamma$ and $\Ckappa$ are introduced,
\begin{align}
    &\varphi_\gamma = \frac 1 2 \gammav \cdot \Cgamma \cdot \gammav , &
    &\varphi_\kappa = \frac 1 2 \kappav \cdot \Ckappa \cdot \kappav ,
    \label{eq:strain_energy_densities}
\end{align}
yielding $\Nv = \Cgamma \cdot \gammav$ and $\Mv = \Ckappa \cdot \kappav$ for the stress resultants.
These tensors can be represented using the conventional expressions for the axial stiffness $\EA$, shear stiffnesses $\GAy$, $\GAz$, the torsional stiffness $\GIt$, and bending stiffnesses $\EIy$, $\EIz$ via
\begin{align}
    \Cgamma &= \EA\, \ev_1 \otimes \ev_1 + \GAy \ev_2 \otimes \ev_2 + \GAz \ev_3 \otimes \ev_3, & 
    \Ckappa &= \GIt\, \ev_1 \otimes \ev_1 + \EIy \ev_2 \otimes \ev_2 + \EIz \ev_3 \otimes \ev_3.
    \label{eq:C_vs_EAI}
\end{align}
This results in the well-known relations for stress $\nv = \rot \cdot \Nv$ and moment stress $\mv = \rot \cdot \Mv$, 
\begin{align}
    \nv &= \left(\EA\, \gamma_1 \right) \gv_1 + \left(\GAy \gamma_2 \right) \gv_2 + \left(\GAz \gamma_3 \right) \gv_3, &
    \mv &= \left(\GIt \kappa_1 \right) \gv_1 + \left(\EIy \kappa_2 \right) \gv_2 + \left(\EIz \kappa_3 \right) \gv_3.
\end{align}
The strain energy stored in the current state of deformation then follows upon integration as
\begin{equation}
    \Phi = \int_{\domain} \left( \varphi_\gamma + \varphi_\kappa \right) ds ,
\end{equation}
and its variation is supposed to equal the virtual work of the internal forces, i.e., $\delta \Phi = \delta \Wint$.

\section{Mixed finite-element formulation}
\label{sec:fe}

This section is dedicated to the design of the proposed beam element. 
The underlying finite element mesh is defined in terms of the arclength coordinate $s$. 
For the moment, we assume a strictly increasing set of nodal coordinates $s_0=0 < s_1 < \dots < s_N = L$, i.e., the beam's centerline does not show any branches. 
Let $\mathcal E = \{E_i = (s_{i-1}, s_i): 1 \leq i \leq N\}$ denote the set of elements.

In their pioneering work, Simo and Vu-Quoc \cite{Simo1986b} propose a finite element formulation based on the consistent discretization of displacements and \emph{incremental} rotations using continuous, piecewise linear shape functions. The element proves capable of reproducing analytic solutions and predicting buckling and post-buckling behavior in a wide range of examples, however, it lacks objectivity. Crisfield and Jeleni\'c \cite{Crisfield1999} identified rotation interpolation in the finite element ansatz as the cause of non-objectivity, and provide a remedy strategy based on a local split of the rotation vector into a constant and superimposed non-constant part. 
For such a split, ensuring continuity of the rotation vector at mesh nodes is not straightforward.
One could, for instance, select the rotation vector of some specific node as constant reference rotation, see, e.g.,~\cite{Crisfield1999,Sonneville2014}, which, however, introduces a mesh-dependency in the interpolation scheme.
In~\cite{Sonneville2017}, an implicit interpolation scheme, similar to that used in geodesic shell elements~\cite{Sander2016}, was employed.


In the present paper, we choose a different strategy. A mixed method, where the moment vector $\Mv$ is treated as independent unknown, is devised. In the context of shell theories, such mixed methods have proven superior in the avoidance of locking \cite{neunteufel2019,neunteufel2021}. While locking compensation may not be the primary motivation for our beam formulation, another trait of mixed methods is intriguing with respect to the design of an objective beam element: in a carefully designed mixed approach, the necessity of continuous approximation is shifted from the rotation vector to the moment vector. Thus, the rotation vector is defined element-local, without need for continuity at mesh nodes. Its continuity is enforced in weak sense through the (extended) principle of virtual work. While our element is of arbitrary polynomial order, we find that the lowest order element sports an element-constant approximation of the rotation vector, with no need for interpolation at all.


To introduce the underlying mixed variational formulation, consider the work of internal forces, as defined through \eqref{eq:delta_W_int} and \eqref{eq:strain_energy_densities},
\begin{align}
    \delta W_{int} = \int_0^L \left( \delta \varphi_\gamma + \delta \varphi_\kappa \right) ds.
\end{align}
To obtain a mixed method, the moment vector $\Mv$ is introduced as an independent field, and the Legendre transform $\overline{\varphi}_\kappa$ of the energy density $\varphi_\kappa$ is needed. In case of a linear material relation \eqref{eq:strain_energy_densities}, this transformation can be computed explicitly as
\begin{align}
    \overline{\varphi}_\kappa 
    = -\frac{1}{2} \Mv \cdot (\mathbb{C}_\kappa)^{-1} \cdot \Mv + \kappav \cdot \Mv .
\end{align}
The variational formulation is now valid for admissible variations of displacement and rotation, but also of the moment vector,
\begin{align}
    \int_0^L \left( \delta \varphi_\gamma + \delta \overline{\varphi}_\kappa \right)\, ds
        = \int_0^L \left( \delta \varphi_\gamma + \Mv \cdot \delta \kappav + (\kappav - (\Ckappa)^{-1} \cdot \Mv) \cdot \delta \Mv\right)  ds.
\end{align}
Above, rotations and their increments are not required to satisfy boundary conditions, while external moments at the beam's tip are considered as essential boundary condition. The correct integration of boundary condition is discussed in detail throughout the remainder of this section.

\subsection{Discrete curvature}
For the definition of curvature $\kappav$ along \eqref{eq:define_kappa}, differentiability of the rotation vectors $\psiv$ and $\psiv_0$ or, equivalently, the rotation tensors $\rot$ and $\rot_0$ with respect to the axial coordinate $s$ is necessary. Unless $\rot$ and $\rot_0$ are not at least (weakly) differentiable, the work integral
\begin{align}
    \int_0^L \kappav \cdot \Mv\, ds \label{eq:work_moment}
\end{align}
is not well-defined in weak sense.
From a geometric point of view, this precludes beams with kinks; in a finite element context, the rotation vector needs to be discretized using continuous, piecewise smooth shape functions. 
In the following, however, a finite element formulation employing discontinuous rotation vector fields shall be developed. Such a formulation directly includes aforementioned beams with kinks, where the beam axis in undeformed configuration is not necessarily smooth, and $\rot_0$ is defined piecewise only. 
To understand the implications of discontinuities in the rotation vector from a mathematical point of view, we shortly explore the notion of \emph{discrete curvature}, which is detailed in the context of shell formulations in \cite{grinspun2006,neunteufel2019}. 

To begin with, we assume the case of a beam with piecewise straight undeformed configuration, where a kink may be present at $s=s_1$. 
We analyze necessary conditions on the moment vector $\Mv$ ensuring a meaningful interpretation of the work integral \eqref{eq:work_moment} for discontinuous, piecewise constant rotation vector fields $\psiv_0$ and $\psiv$, 
\begin{align}
    \psiv_0(s) &= \left\{ \begin{array}{ll} \psiv_{0,1}&\text{if } s < s_1,\\ \psiv_{0,2}& \text{if } s > s_1,\end{array} \right. &
    \psiv(s) &= \left\{ \begin{array}{ll} \psiv_1&\text{if } s < s_1,\\ \psiv_2& \text{if } s > s_1.\end{array} \right.
\end{align}
To this end, we define a family of smooth rotation vector fields $\psiv_{0,\epsilon}$ and $\psiv_\epsilon$ that converge to $\psiv_0$ and $\psiv$, respectively, as $\epsilon \to 0$ in the following. The limiting value of the work integral shall then motivate the definition of discrete curvature. 

We find the relative rotation $\psiv_{12}$ between the two beam sections via
\begin{align}
    \psiv_{12} := \sphericalangle (\rotmat_1^T \cdot \rotmat_2),
\end{align}
where we used the shorthand notation $\rot_i = \exp(\skw\psiv_i)$ for $i=1,2$.
For the undeformed configuration, the relative rotation $\psiv_{0,12}$ is defined analogously.
It is important to note that this relative rotation is objective, i.e., for a constant rigid body rotation $\rot_{rb}$ that is superimposed to both $\bar{\rot}_i = \rot_{rb} \cdot \rot_i$, the relative rotation remains unchanged,
\begin{align}
    \bar{\psiv}_{12} 
    = \sphericalangle (\bar{\rotmat}_1^T \cdot \bar{\rotmat}_2) 
    = \sphericalangle ({\rotmat}_1^T \cdot {\rot}_{rb}^T \cdot {\rot}_{rb} \cdot {\rotmat}_2) 
    = \sphericalangle (\rotmat_1^T \cdot \rotmat_2) 
    = \psiv_{12}.
\end{align}
Using the relative rotation vector, a smooth approximation of the rotation tensor field is given through
\begin{align}
    \rotmat_\epsilon(s) := \left\{
        \begin{array}{ll}
            \rotmat_1 & \text{if } s \leq s_1,\\
            \rotmat_1 \cdot \exp(\frac{s-s_1}{\epsilon} \skw{\psiv}_{12}) & \text{if } s_1 < s \leq s_1+\epsilon,\\
            \rotmat_2 & \text{if } s > s_1+\epsilon.
        \end{array} 
    \right.
\end{align}
%
By definition, $\rot_\epsilon$ is continuous at $s=s_1$ and $s = s_1+\epsilon$, and otherwise smooth. We use $\psiv_\epsilon = \sphericalangle(\rotmat_\epsilon)$ as smooth approximation to the discontinuous rotation vector field $\psiv$; note that, in general, $\psiv_\epsilon \neq \psiv_1 + \left(\psiv_2 - \psiv_1 \right) (s-s_1)/\epsilon $. An analogous procedure yields the smooth rotation vector field $\psiv_{0,\epsilon}$ for the undeformed configuration.

The vector of moment strains ${\kappav}_\epsilon$ is well-defined in classical sense, its representation as skew-symmetric tensor is,
\begin{align}
    \skw\kappav_\epsilon = \rotmat_\epsilon^T \cdot \rotmat_\epsilon' - \rotmat_{0,\epsilon}^T \cdot \rotmat_{0,\epsilon}'.
\end{align}
It depends only on the relative angles $\psiv_{12}$ and $\psiv_{0,12}$, as can easily be verified: outside $(s_1, s_1+\epsilon)$, $\kappav_\epsilon$ vanishes as $\rotmat_\epsilon$ and $\rotmat_{0,\epsilon}$ are constant. For $s \in (s_1, s_1+\epsilon)$, we find that, since the direction of the axis of relative rotation $\psiv_{12}$ is constant throughout the interval, algebraic manipulation yields
\begin{align}
    \rotmat_\epsilon^T \cdot \rotmat_\epsilon' 
    = \exp\left(\frac{s-s_1}{\epsilon} \skw{\psiv}_{12}\right)^T \cdot \frac{d}{ds}\exp\left(\frac{s-s_1}{\epsilon} \skw{\psiv}_{12}\right)
    = \frac{1}{\epsilon} \skw\psiv_{12}.
\end{align}
A similar representation can be found for the initial configuration based on $\rotmat_{0,\epsilon}$ and $\psiv_{0,12}$.
The moment strain vector finally is given by
\begin{align}
    \kappav_\epsilon &= \frac{1}{\epsilon} \left(\psiv_{12} - \psiv_{0,12}\right) \qquad \text{for } s \in (s_1, s_1+\epsilon).
\end{align}
For some moment vector field $\Mv$, the work integral computes as
\begin{align}
    \int_0^L \kappav_\epsilon \cdot \Mv\, ds 
                       = \frac{1}{\epsilon}\left(\psiv_{12} - \psiv_{0,12}\right) \cdot \int_{s_1}^{s_1+\epsilon} \Mv\, ds.
\end{align}
For a \emph{continuous} moment vector $\Mv$, the mean value theorem guarantees the existence of some $\bar s_\epsilon \in [s_1, s_1+\epsilon]$ such that
\begin{align}
    \int_0^L \kappav_\epsilon \cdot \Mv\, ds = \left(\psiv_{12} - \psiv_{0,12}\right) \cdot \Mv|_{\bar s_\epsilon}.
\end{align}
As $\epsilon \to 0$, the intermediate position $\bar s_\epsilon \to s_1$, and we find the representation of the work integral as the inner product of relative rotation and moment vector,
\begin{align}
    \int_0^L \kappav_\epsilon \cdot \Mv\, ds \to \left(\psiv_{12} - \psiv_{0,12}\right) \cdot \Mv|_{s_1}. \label{eq:work_moment_disc}
\end{align}
Thus, when using the concept of discrete curvature, local moments are dual to the \emph{change} of relative rotation observed in the transition from one side of the interface to the other. If $\psiv_{0,12} \neq \boldsymbol{0}$, i.e., if there is a kink at position $s_1$ in the undeformed configuration, this angle must be retained in the deformed configuration for a vanishing discrete moment strain.
Note that a consistent definition of the moment vector at interface point $s_1$ is necessary. In some sense, the usage of discrete curvature and moment strains transfers the necessity of inter-element continuity from the rotation vector to the moment vector.

The above considerations can be extended to piecewise smooth but discontinuous rotation vector fields $\psiv$ and $\psiv_0$, which are relevant in a finite element discretization. 
Let $\Mv$ denote the moment vector, which is continuous across all element interfaces, i.e., in all mesh nodes $s_i$, and let $\rot_i = \exp(\skw\psiv)|_{E_i}$ denote the rotation tensor restricted to element $E_i$. Using $\sphericalangle(\rot_{i-1}^T \cdot \rot_i)$ and $\sphericalangle(\rot_{0,i-1}^T \cdot \rot_{0,i})$ accordingly for the relative rotation across element interfaces, \eqref{eq:work_moment} can be replaced by
\begin{align}
    \sum_{i=1}^{N} \int_{E_i} \kappav \cdot \Mv\, ds + \sum_{i=0}^{N} \left(\sphericalangle(\rot_{i-1}^T \cdot \rot_i) - \sphericalangle(\rot_{0,i-1}^T \cdot \rot_{0,i})\right) \cdot \Mv|_{s_i}. \label{eq:work_moment_fe_disc}
\end{align}
Above, the beam's tip points at $s_0=0$ and $s_N=L$ have been included tacitly; at any boundary vertex, boundary conditions are included setting either $\Mv$, or, alternatively, defining the relative rotation $\sphericalangle(\rot_{i-1}^T \cdot \rot_i)$ as relative to the prescribed boundary value.

\subsection{A mixed beam element}


Using the above notion of discrete curvature, we can define a finite-element method that uses continuous displacements and moments along with discontinuous rotations as degrees of freedom. 
In such a method, however, kinematic boundary conditions on rotations, which are required to realize a clamped end, for instance, cannot be included explicitly. 
Instead, they need to be imposed in weak form through the moment vector, which is unconstrained for boundary conditions of this kind.
From an implementational point of view, though, the weak representation of kinematic boundary conditions does not present any difficulties.
External moments, on the other hand, need to be included explicitly as essential boundary conditions, which is certainly challenging due to the fact that material moments $\Mv$ are interpolated rather than actual moments $\mv = \rot \cdot \Mv$, cf.~\eqref{eq:material_stress_resultants}.
For this reason, we introduce nodal rotations, which are only defined at the nodes of the mesh, as an additional unknown field.

On a given finite element mesh, the beam element uses the following unknown fields:
\begin{itemize}
    \item positions $\rv$ (or displacements $\uv = \rv - \rv_0$) interpolated by continuous, piecewise polynomial functions of order~$k$,
    \item element-local moments $\Mv$ represented by discontinuous, piecewise polynomial functions of order~$k$,
    \item nodal \emph{incremental} rotations $\psiv^{(V)}_i$, single-valued in all mesh nodes, with no interpolation in the element interior, 
    \item element-local rotations $\psiv^{(E)}$ represented by discontinuous, piecewise polynomial functions of order~$k-1$.
\end{itemize}
On a single element $E_i = (s_{i-1}, s_i)$, the work of internal forces is computed both through integration over the axial coordinate, as well as through element-boundary terms,
\begin{align}
    \int_{E_i} \left( \varphi_\gamma + {\overline{\varphi}}_\kappa \right)\, ds
    + \left[\sphericalangle((\rot^{(E)})^T \cdot \rot^{(V)}) \cdot \Mv\right]_{s_{i-1}}^{s_i}.
\end{align}
Above and in the following, we denote the rotation tensor computed from the element-local rotation vector as $\rot^{(E)} := \exp(\skw\psiv\vphantom{\psi}^{(E)})$, while the rotation tensor due to nodal incremental rotations is defined in both element nodes as $\rot^{(V)} := \exp(\skw{\psiv}\vphantom{\psiv}^{(V)}) \cdot \exp(\skw{\psiv}_{0})$. 
We emphasize that in the computation of (variations of) the potentials $\varphi_\gamma$ and $\varphi_\kappa$, only the element-local rotations $\psi^{(E)}$ are used, such that no interpolation of the nodal incremental rotations $\psiv^{(V)}$ to the element's interior is necessary.

The introduction of independent $\psiv^{(V)}$ can be interpreted as a \emph{hybridization strategy}, where the nodal rotations act as a Lagrangian multiplier enforcing the broken continuity of the moment vector. A similar procedure is proposed in the context of shells by Neunteufel and Sch\"oberl \cite{neunteufel2019,neunteufel2021}. From a computational point of view, hybridization is often advantageous as it facilitates static condensation, which, in turn, leads to a positive definite system matrix of reduced size. For the proposed beam element, only displacements and nodal incremental rotations remain as non-condensable degrees of freedom in the global system matrix. 
While the increase in computational efficiency may be significant for volumetric elements, beam elements are connected solely via their nodes, so static condensation plays only a minor role in this respect. 
The main advantage of the hybridization is due to the fact that nodal incremental rotations allow for a straightforward representation of classical boundary conditions: for clamped conditions, requiring $\psiv^{(V)} = \zerov$ explicitly is sufficient; when prescribing rotations, this can be done directly for the nodal values. 
As far as external couples are concerned, the virtual work of external forces can be evaluated using both element-local and nodal rotations,
\begin{equation}
    \delta \Wext = \int_0^L \left( \qv \cdot \delta \rv + \cv \cdot \delta \thetav^{(E)} \right) ds
    + \left[ \nv \cdot \delta \rv + \mv \cdot \delta \thetav^{(V)} \right]_0^L ,
\end{equation}

In the context of shells, this hybridization procedure of introducing $\psiv^{(V)}$ ensures that the correct balance of moments is enforced also along edges where several branches of the structure meet \cite{neunteufel2019}. The same benefit is observed for the present beam element: T-joints or branching points are modeled correctly without further treatment. The constraint induced by $\psiv^{(V)}$ in a single node is equivalent to the balance of moments around that node, regarding contributions from all branches emanating there equally.

To ensure objectivity of moment strains, which are computed from element-local rotation vectors $\psiv^{(E)}$, we introduce a multiplicative splitting of the corresponding rotation tensors $\rot^{(E)} := \exp(\skw{\psiv}\vphantom{\psiv}^{(E)})$.
In line with the reasoning of Crisfield and Jelenić~\cite{Crisfield1999}, total rotations are represented as an element-wise constant rotation (the lowest-order contribution) superimposed by relative rotations with zero mean value (the higher-order contribution). Let $\rot_{lo}^{(E)} := \exp (\skw{\rotvec}\vphantom{\psiv}_{lo}^{(E)})$ and $\rot_{ho}^{(E)} := \exp( \skw{\rotvec}\vphantom{\psiv}_{ho}^{(E)} )$ denote the respective contributions to the rotation tensors; the total rotation is given by the successive application
\begin{equation}
    \rot^{(E)} = \rot_{lo}^{(E)} \cdot \rot_{ho}^{(E)} .
    \label{eq:rot_lo_ho}
\end{equation}
The constant rotation is represented by element-wise constant rotation vectors $\rotvec_{lo}^{(E)}$; for the higher-order relative rotations $\rotvec_{ho}^{(E)}$, we use piecewise polynomial functions of order $k-1$ with zero mean value. This split is natural when using Legendre polynomials as a basis, where the first shape function corresponds to the constant field, while all remaining shape functions have zero mean value due to orthogonality. 
Note that for the lowest-order approximation $k=1$, the element-local rotation equals the constant part, i.e., $\rot^{(E)} = \rot_{lo}^{(E)}$.

\subsection{Locking behavior}
\label{sec:locking}
Shear and membrane locking are effects that are often observed to increase with decreasing thickness $t$ of the beam. As the ratio of cross-section diameter and element length decreases, axial and shear stiffnesses collected in $\Cgamma$ dominate bending and torsional stiffnesses in $\Ckappa$ by $t^{-2}$. In the limit case of an inextensible Bernoulli-Euler beam of vanishing thickness, the constraint $\gammav = \boldsymbol{0}$ is enforced; small but positive thickness is often compared to a penalty realization of this constraint. \emph{Locking} occurs whenever the discretization of different quantities does not allow for the constraint $\gammav = \boldsymbol{0}$ to be satisfied exactly, such that spurious strains arise that lead to parasitic energy contributions, as is detailed, e.g., by Koschnick et al.~\cite{koschnick2005}. 

We motivate whether locking is expected for the proposed element. The definition of the force strain~\eqref{eq:gammav}
can be satisfied exactly if and only if rotation $\psiv$ and deformation $\rv$ match in such a way that constant $\rot^T  \cdot \rv' = \ev_1$ can be achieved, or equivalently $\rv' = \rot \cdot \ev_1 = \gv_1$. This is clearly the case for the lowest-order method with $k=1$: then, $\rv'$ is element-wise constant, and element-wise constant $\psiv^{(E)}$ can be chosen to meet the constraint exactly.

For higher-order methods, $k\geq 2$, this relation no longer holds due to the nonlinear relation of rotation $\rot$, rotation vector $\psiv$, and its representation via $\psiv^{(E)}_{lo}$ and $\psiv^{(E)}_{ho}$. To preserve a locking-free method of optimal order, a projection of the stress resultant is often used. In their early work, Stolarski et al.~\cite{stolarskl1983} showed that reduced integration in the shear term is appropriate for shear-deformable beams, provided the decomposition of shear and bending contributions is exact. For the proposed beam element, we choose a similar approach of reduced integration in the strain energy $\varphi_\gamma$. When using Gauss integration points, this is equivalent to projecting force strains $\gammav$ to $P^{k-1}$ locally. Since displacements are discretized at order $k$, we see that also $\rv'$ is in $P^{k-1}$ locally, so the order of approximation is not diminished by reduced integration. 

A numerical assessment of the locking behavior and the effect of reduced integration in different contributions to the virtual work statement is provided in Sections~\ref{sec:examples_rot2circle} and~\ref{sec:arc_Bathe}.

\section{Numerical examples}
\label{sec:examples}
In what follows, we present several numerical examples which are meant to show the potential of the proposed beam element. 
Except for the last one, all problems considered subsequently are typical benchmarks in the context of large-deformation beam theories. Consistent units are assumed for all dimensional quantities.
All numerical examples are implemented in the general-purpose finite-element code \texttt{Netgen/NGSolve};\footnote{\href{https://www.ngsolve.org}{\texttt{https://www.ngsolve.org}}} implementations of our examples are available in a public repository.\footnote{\href{https://gitlab.com/alexander.humer/geometrically\_exact\_beam}{\texttt{https://gitlab.com/alexander.humer/geometrically\_exact\_beam}}}

\subsection{Rod bent to circle}
\label{sec:examples_rot2circle}
\begin{figure}[htb!]
    \centering
    \includegraphics{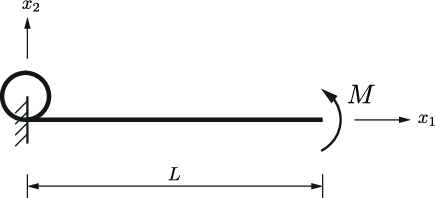}
    \caption{Roll-up of a cantilever beam under a tip-moment: Undeformed, straight configuration and deformed shape of a double circle corresponding to $M = 4 \pi EI / L$.}
    \label{fig:1_circle}
\end{figure}

As a first example, we study the planar problem of a cantilever beam subject to a tip moment such that it bends into a (double) circle, see Fig.~\ref{fig:1_circle}.
The problem was already analyzed by Simo and Vu-Quoc~\cite{Simo1986b} in their pioneering work, which can be considered the first successful finite-element formulation for geometrically exact beams.
Ever since, the problem has become one of the benchmarks novel approaches need to stand up to, see, e.g.,~\cite{Sonneville2014, Meier2015, harsch2023, Wasmer2024} for recent contributions considering different variants.
The purpose of the problem is to study how well a pure state of bending can be represented by the proposed beam element using different polynomial orders.
We also study how locking phenomena can be alleviated by means of reduced integration in selected contributions to the strain energy.

To quantify the accuracy of approximate solutions, we use the normalized (by the beam's length) $L^2$-error in the centerline position $\rv$ with respect to the reference solution $\rv_{\mathrm{ref}}$, 
\begin{equation}
    \err = \frac{1}{L} \sqrt{\frac 1 L \int_0^L \left\Vert \rv - \rv_{\mathrm{ref}}\right\Vert^2 ds} ,
\end{equation}
which is given in closed-form for the present problem.

For linear constitutive relations~\eqref{eq:constitutive_equations}, planar problems of shear-deformable beams are governed by two non-dimensional similarity parameters, which relate the beam's axial, shear, and bending stiffnesses.
For the present problem, we assume a homogeneous, square cross-section of side length $a$, which defines the slenderness of the beam, i.e., the ratio between axial and lateral dimensions, as $\rho := L / a$.
In terms of the slenderness, the non-dimensional ratio between axial stiffness and bending stiffness is then given by $L^2 \EA / \EIz = 12 \rho^2$.
To bend the initially straight beam into a double circle, a tip moment $\Mv = M \ez$, with $M = 4 \pi \EIz / L$, is applied at the free end ($s = L$).

To study how well the proposed formulation is able to represent a pure state of bending, we first vary the slenderness from (moderately) thick ($\rho = 10$) to slender beams ($\rho = 1000$) and study the evolution of the error as the number of elements increases ($\nel = 8 \ldots 512$, in powers of $2$).
We note that as few as $4$ elements per circle are used in the coarsest setting. 
The polynomial order of elements is varied from the lowest-order case ($k = 1$) up to fourth-order elements ($k=4$).

\begin{figure}[htb!]
    \centering
    \includegraphics[width=1.\textwidth]{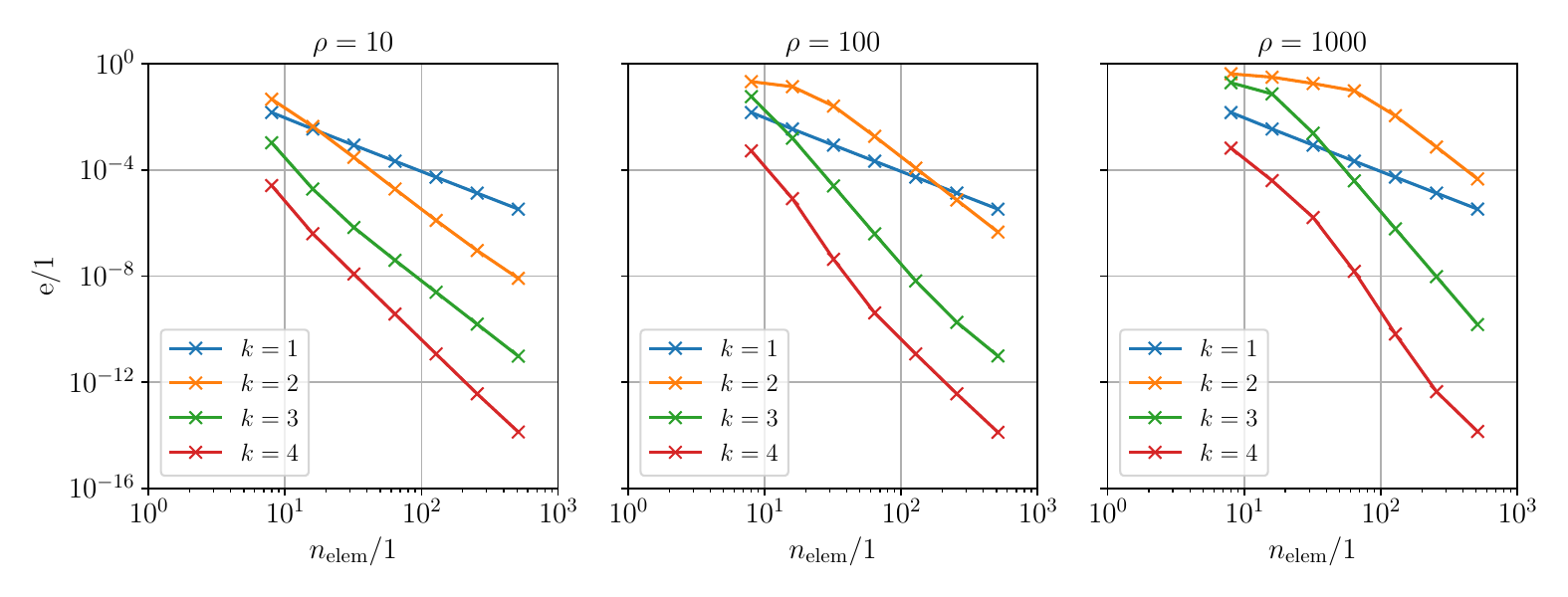}
    \caption{Roll-up of a cantilever beam under a tip-moment: rates of convergence for beam elements of orders $k=1,\ldots,4$ and different values of slenderness; $\rho = 10$ (left), $\rho = 100$ (center), and $\rho = 1000$ (right).}
    \label{fig:1_circle_rho}
\end{figure}

For thick beams with a slenderness of $\rho = 10$ (left), we obtain close to optimal convergence rates in the $L^2$-error, i.e., $O(h^{k+1})$ for elements of order $k$ and element size $h = 1/n_\text{elem}$.
Interestingly, lowest-order elements are more accurate than the quadratic elements for very coarse meshes.
Once we increase the slenderness, however, the accuracy starts to deteriorate. 
For a slenderness of $\rho = 100$ (center), the convergence of quadratic elements already suffers significantly; for very slender beams characterized by $\rho = 1000$, this gets even more pronounced, and also cubic elements are visibly affected, see Fig.~\ref{fig:1_circle_rho} (right).
Fourth-order elements remain largely unaffected by the slenderness.
Only for coarse meshes, we can observe an increase in the error.
Remarkably, the error obtained for the lowest-order approximation is completely independent of the slenderness. These results are in accordance with the findings of Section~\ref{sec:locking}, where we argued that locking is expected only for higher-order discretizations.
%

\begin{figure}[htb!]
    \centering
    \includegraphics[width=1.\textwidth]{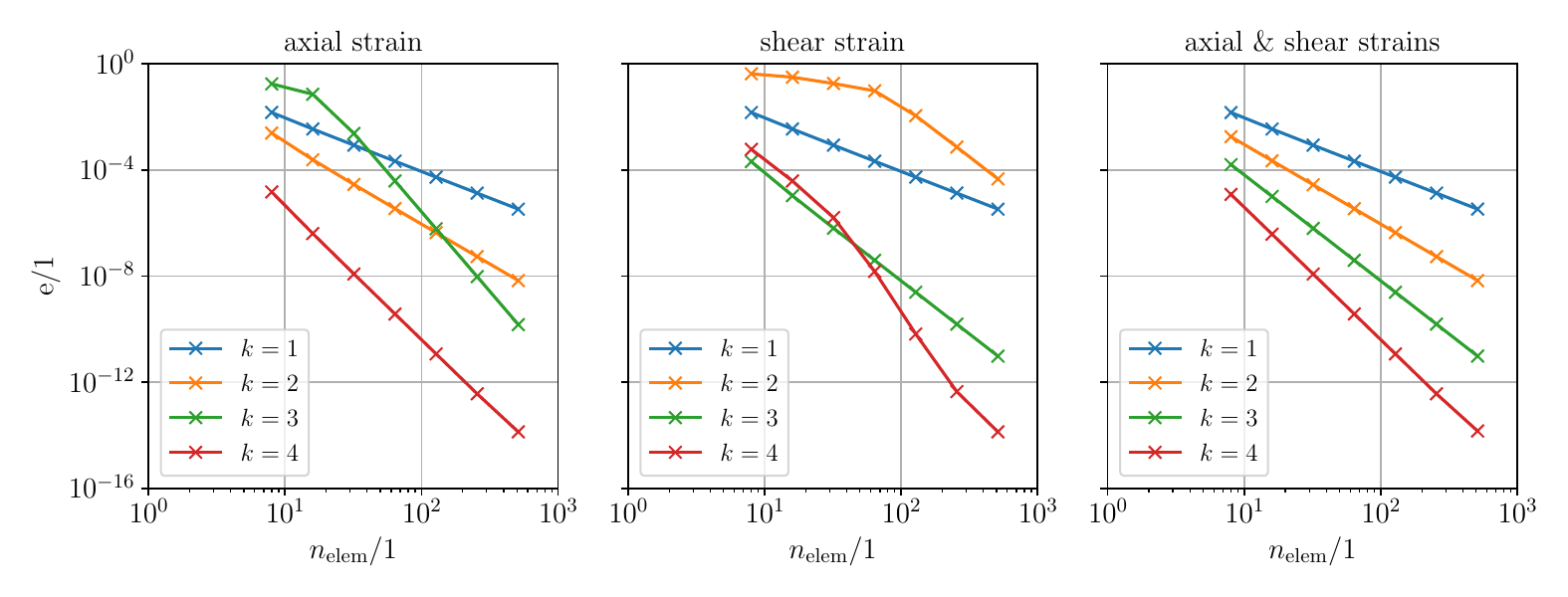}
    \caption{Roll-up of a cantilever beam under a tip-moment: rates of convergence when applying reduced integration in the strain energy related to axial strains (left), shear strains (center), and all force strains (right). The slenderness is set to $\rho = 1000$.}
    \label{fig:1_circle_red}
\end{figure}

Except for the lowest-order element, a spurious coupling between moment and force strains is present in a state of pure bending, which constitutes what is commonly referred to as \emph{locking}.
To further analyze the specific nature of locking, we selectively apply reduced integration, which has proven an effective and simple means to reduce locking phenomena, to individual terms of the strain energy.
In particular, we successively underintegrate---by one order---the strain energy related to extension first, then only the contributions due to shear deformation, and eventually both terms related to axial and to shear deformation, see Fig.~\ref{fig:1_circle_red}.
We stay with the case of $\rho = 1000$, since locking becomes more severe as the slenderness increases.

Reduced integration in the strain energy associated with axial force strains improves the accuracy of second- and fourth-order elements, see Fig.~\ref{fig:1_circle_red} (left).
Indeed, we recover ideal convergence rates by reducing the coupling between axial extension and bending deformation, which suggests that elements of even order suffer from \emph{membrane locking}. 
For cubic elements, no improvements are observed.

Next, we underintegrate the contribution to strain energy related to shear deformation, which entirely reverses the picture, see Fig.~\ref{fig:1_circle_red} (center).
Now, the cubic element shows an ideal rate of convergence, whereas second- and fourth-order elements revert to the results obtained with full integration, which indicates that the former suffers from \emph{shear locking}.

Reduced integration in all terms of the strain energy $\varphi_\gamma$, cf. Eq.~\eqref{eq:strain_energy_densities}, effectively eliminates all locking phenomena for the element orders considered herein. 
Once again, we emphasize that the lowest-order element ($k=1$) is a notable exception since it does not show any locking irrespective of the integration. 


\subsection{Superimposed rotation of deformed cantilever}

\begin{figure}[htb!]
    \centering
    \includegraphics[height=4cm]{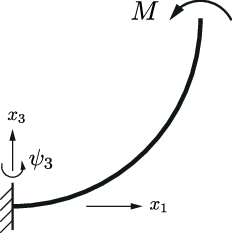}
    \hspace{1cm}
    \includegraphics[height=4.5cm,clip=true, trim=0 0.5cm 0 0.5cm]{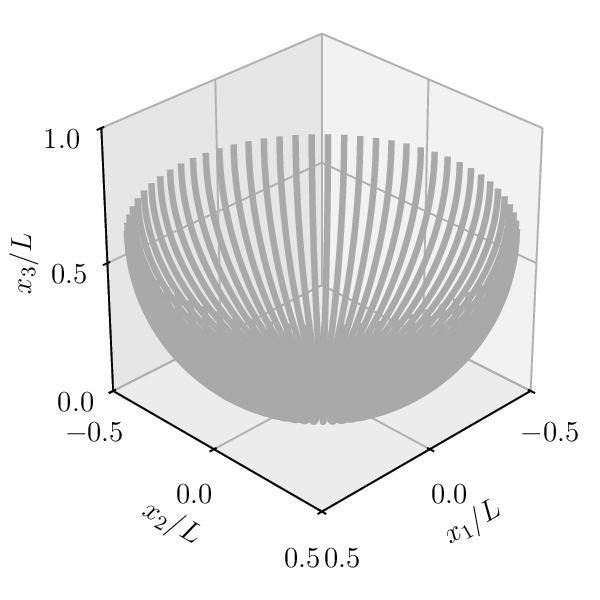}
    \caption{Superimposed rotation of deformed cantilever: deformed configuration prior to the rigid-body rotation (left); snapshots of deformed configurations in the course of ten rotations (right).}
    \label{fig:2_objectivity}
\end{figure}

In the proposed formulation, generalized (moment) strains are computed from the high-order part~$\rot_{ho}^{(E)}$ of element-local rotations, cf. Eq.~\eqref{eq:rot_lo_ho}.
As argued by Crisfield and Jelenić~\cite{Crisfield1999}, any rigid-body motion superimposed on the deformed configuration would not affect the moment strains also in the discrete setting.  
The objectivity of the proposed beam element can be verified by a simple numerical test, which was also used in several previous contributions, see, e.g.,~\cite{Meier2019,harsch2023}. 
We replicate the defining feature that constitutes objectivity by first bending a cantilever beam into a quarter circle using a tip-moment $M = -\pi \EIy / 2 L$ about the $x_2$-axis, before the deformed beam is subjected to a rigid-body rotation about the $x_3$-axis, see Fig.~\ref{fig:2_objectivity} (left).
To impose the rotation, the nodal rotation at the clamped end of the cantilever $\psiv^{(V)}(0)$ is prescribed kinematically. 
For the present study, we use the same properties as in the previous problem (Section~\ref{sec:examples_rot2circle}); the slenderness is set to $\rho = 10$.

Under the superimposed rigid-body rotation, the transverse tip-displacement $u_3(L)$ is expected to remain constant, and so is the stored elastic energy, which, for the present problem of an initial bending about the $x_2$-axis, is given by
\begin{equation}
    \overline{\Phi}_\kappa = \int_0^L \frac{(M_2)^2}{\EIy} ds .
\end{equation}

We study ten full rotations about the $x_3$-axis, which we impose in \num{71} increments for the sake of visualization, see Fig.~\ref{fig:2_objectivity} (right) for the respective deformed configurations obtained with four cubic elements.

\begin{figure}[htb!]
    \centering
    \includegraphics[width=1.\textwidth]{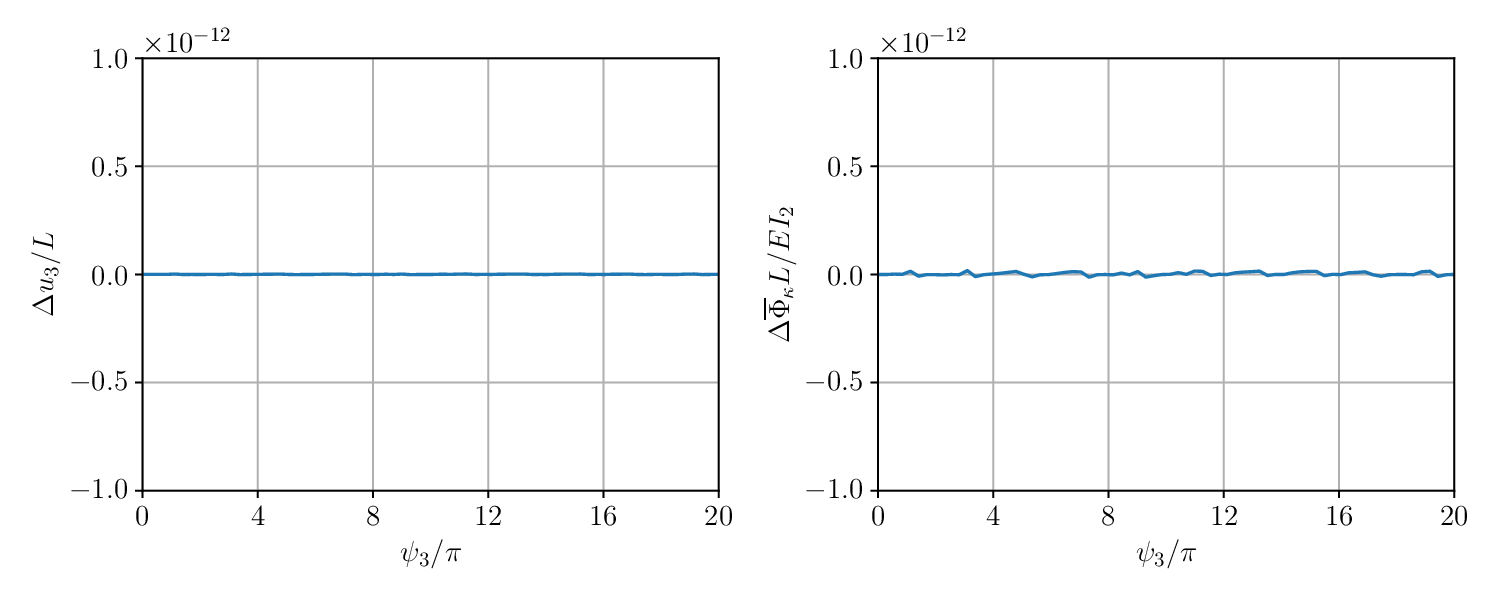}
    \caption{Superimposed rotation of deformed cantilever: Deviation of the transverse tip-displacement and the potential energy from the initial state as a function of the prescribed angle of rotation.}
    \label{fig:2_objectivity_res}
\end{figure}

To quantify a potential loss of objectivity, we consider the deviation of both the transverse tip-displacement $u_3(L)$ and the potential energy $\overline{\Phi}_\kappa$ from the bent state prior to the rotation, which are denoted by $\Delta u_3$ and $\Delta \overline{\Phi}_\kappa$ subsequently. 
Figure~\ref{fig:2_objectivity_res} illustrates $\Delta u_3$ and $\Delta \overline{\Phi}_\kappa$ as functions of the prescribed rotation $\psi_3$.
Clearly, both quantities remain preserved up to floating-point (double) precision with the convergence criterion of Newton's method set accordingly.
We can therefore conclude that the proposed beam element satisfies the requirement of objectivity.

\subsection{Arc-segment under out-of-plane force}
\label{sec:arc_Bathe}

\begin{figure}[htb!]
    \centering
    \includegraphics[width=0.3\linewidth,trim=0 0 0 1mm,clip=true]{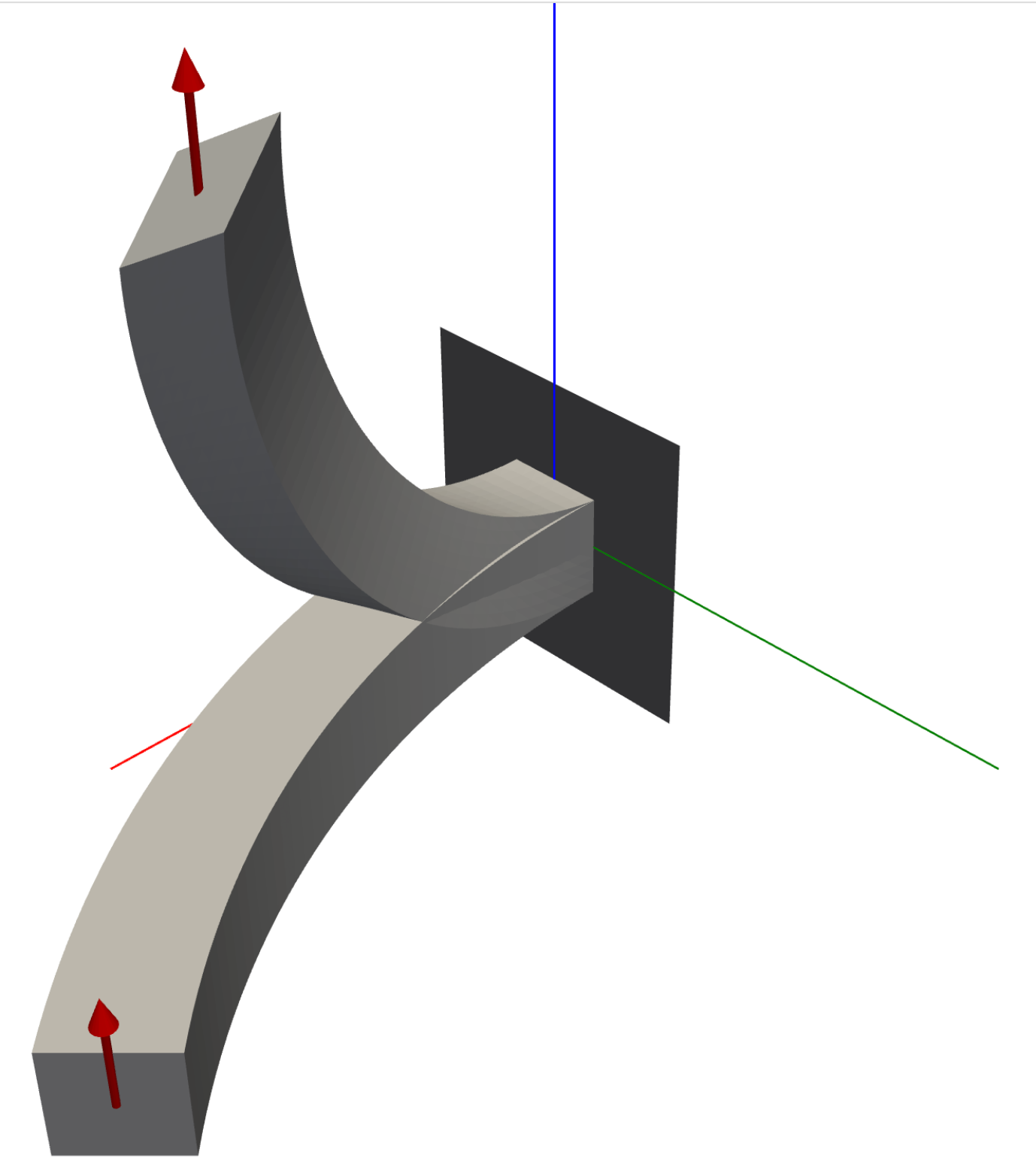}
    \includegraphics[width=0.3\linewidth,trim=0 0 0 1mm,clip=true]{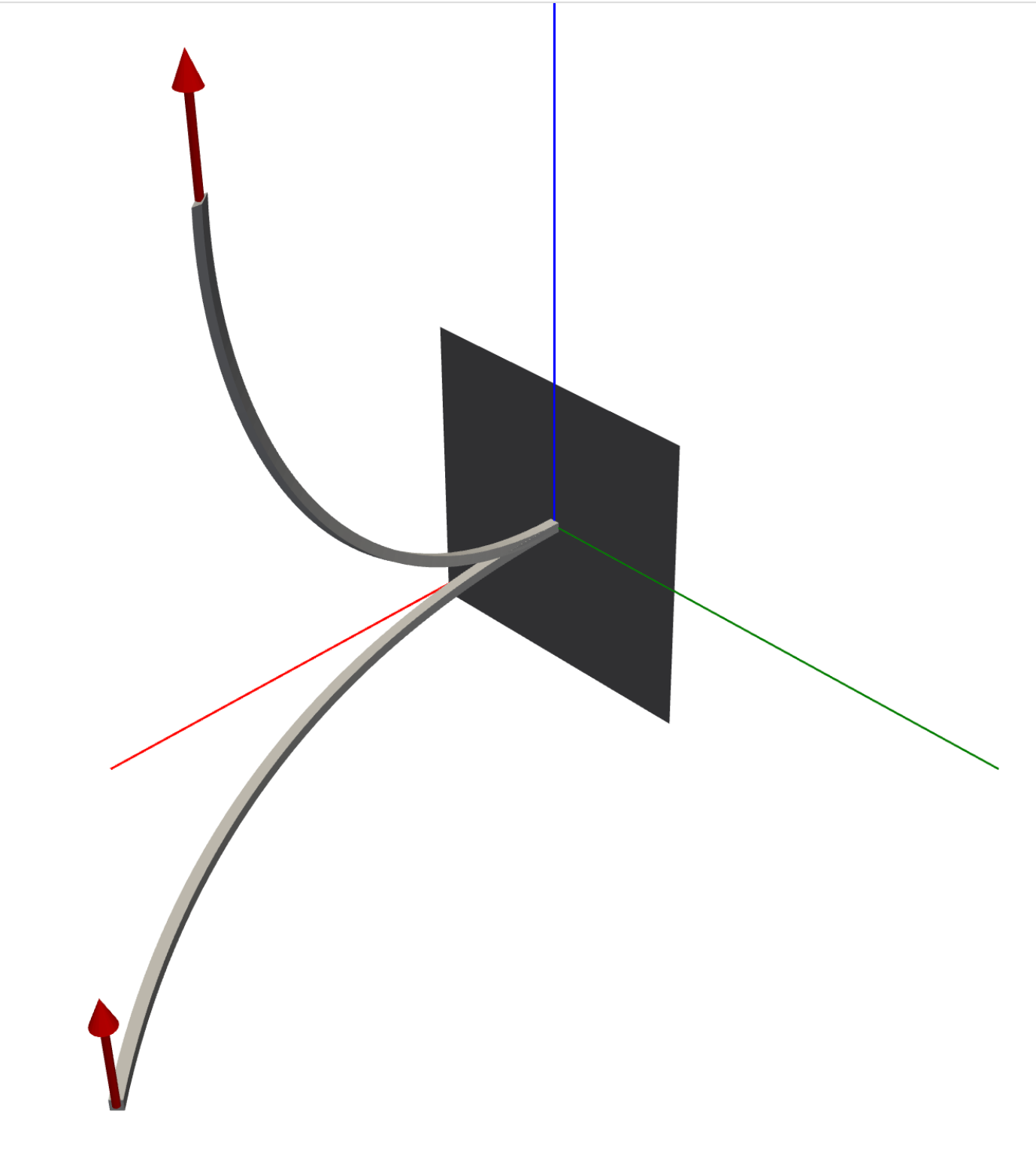}
    \includegraphics[width=0.3\linewidth,trim=0 0 0 1mm,clip=true]{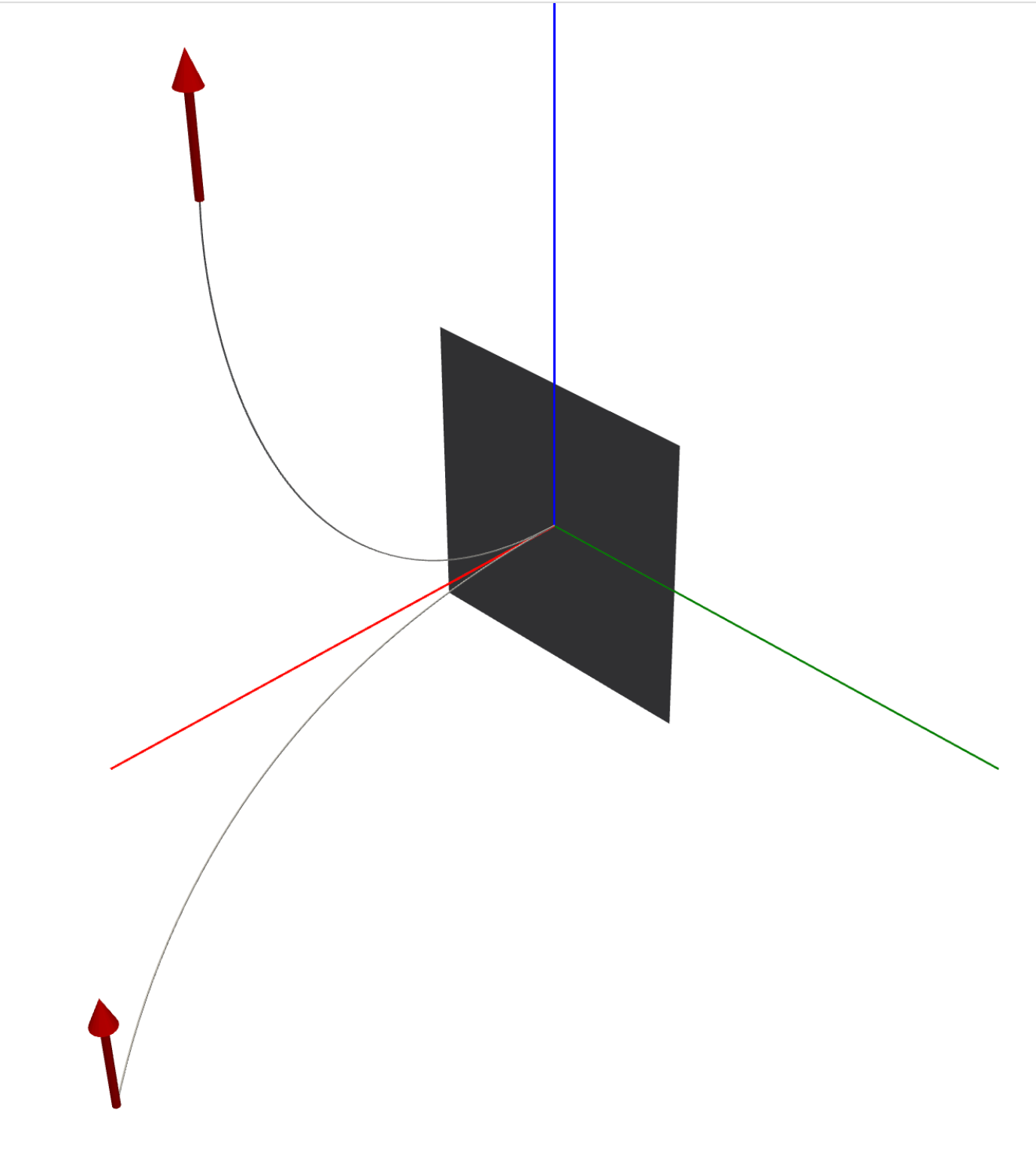}
    \caption{Arc-segment under out-of-plane force: undeformed and deformed configurations for slenderness $\rho=10$ (left), $\rho=100$ (center), and $\rho=1000$ (right).}
    \label{fig:3_arc_configs}
\end{figure}

In the previous examples, we have studied beams subjected only to bending deformation (and rigid-body motion).
Next, we consider another classic benchmark originally proposed in an early finite-element formulation for large deformation problems by Bathe and Bolourchi~\cite{Bathe1979}, where an octant of a circle (\qty{45}{\degree}-bend) is subjected to an out-of-plane force.
The arc with a radius of curvature $R=100$ lies within the $x_1x_2$-plane of the global frame and it is clamped at $s=0$.
At the free end, a concentrated force perpendicular to the arc's plane $\nv(L) = F \ez$ is applied.
Assuming square cross-sections of side length $a$, we use a slightly modified definition of the slenderness in terms of the radius, $\rho = R / a$, for consistency with previous contributions~\cite{Simo1986,Ibrahimbegovic1995,Jelenic1999,Meier2019}.
Again, we study the influence of slenderness, ranging from (moderately) thick ($\rho = 10$) to very slender ($\rho=1000$) beams, see Fig.~\ref{fig:3_arc_configs}.

The effective cross-sectional stiffnesses of the beam are set to $\EA = \num{1e7}a^2$, $\GAy = \GAz = \EA / 2$, and $\GIt = \EIy = \EIz = \num{1e7} a^4/12$.
To maintain a comparable displacement irrespective of the slenderness, the external load is set to $F = \num{600} \left(100/\rho\right)^4$.

To mitigate locking, we employ reduced integration in all terms of the force-strain energy $\varphi_\gamma$, which has proven effective in our first example, cf.~Section~\ref{sec:examples_rot2circle} above.
Figure~\ref{fig:3_arc_configs} shows the undeformed and deformed configurations for three different slenderness values $\rho \in \left\{ 10, 100, 1000\right\}$.
The corresponding numerical values of the tip-displacement $\uv(L)$ are listed in Tab.~\ref{tab:3-arc_results}.

\begin{table}[htb!]
    \sisetup{
    round-mode = figures,
    round-precision = 6,
    group-digits = integer
    }
    \centering
    \caption{Arc-segment under out-of-plane force: components of the tip-displacement $\uv(L)$ for slendernesses $\rho \in \left\{ 10, 100, 1000\right\}$, element orders $k=1,\ldots,4$, and number of elements $\nel \in \left\{ 4, 32\right\}$.}
    \label{tab:3-arc_results}
    \begin{tabular}{c c c c c c c c c c c}
        \toprule
        &   &   \multicolumn{3}{l}{$\rho=10$} 
            &   \multicolumn{3}{l}{$\rho=100$}
            &   \multicolumn{3}{l}{$\rho=1000$}\\
        \cmidrule(lr){3-5}
        \cmidrule(lr){6-8} 
        \cmidrule(lr){9-11} 
        $k$ &   $\nel$  
        &   $u_1(L)$    &   $u_2(L)$    &   $u_3(L)$    
        &   $u_1(L)$    &   $u_2(L)$    &   $u_3(L)$    
        &   $u_1(L)$    &   $u_2(L)$    &   $u_3(L)$ \\
        \midrule
        1   &   4   &   \num{-23.75125} & \num{-13.59944} & \num{54.49189} 
                    &   \num{-23.66304} & \num{-13.57181} & \num{53.87134} 
                    &   \num{-23.66215} & \num{-13.57153} & \num{53.86514} \\
        2   &   4   &   \num{-23.65101} & \num{-13.63756} & \num{54.08983} 
                    &   \num{-23.56625} & \num{-13.61001} & \num{53.46957} 
                    &   \num{-23.56541} & \num{-13.60974} & \num{53.46337} \\
        3   &   4   &   \num{-23.64497} & \num{-13.63205} & \num{54.095} 
                    &   \num{-23.5602} & \num{-13.6045} & \num{53.47483} 
                    &   \num{-23.55935} & \num{-13.60423} & \num{53.46863} \\
        4   &   4   &   \num{-23.64501} & \num{-13.63207} & \num{54.09503} 
                    &   \num{-23.56024} & \num{-13.60452} & \num{53.47486} 
                    &   \num{-23.55939} & \num{-13.60425} & \num{53.46866} \\
        \midrule
        1   &   32  &   \num{-23.64685} & \num{-13.63184} & \num{54.10082} 
                    &   \num{-23.56202} & \num{-13.60429} & \num{53.48064} 
                    &   \num{-23.56118} & \num{-13.60401} & \num{53.47444} \\
        2   &   32  &   \num{-23.64503} & \num{-13.63208} & \num{54.09503} 
                    &   \num{-23.56026} & \num{-13.60453} & \num{53.47486} 
                    &   \num{-23.55941} & \num{-13.60426} & \num{53.46866} \\
        3   &   32  &   \num{-23.64501} & \num{-13.63207} & \num{54.09503} 
                    &   \num{-23.56024} & \num{-13.60452} & \num{53.47486} 
                    &   \num{-23.55939} & \num{-13.60425} & \num{53.46866} \\
        4   &   32  &   \num{-23.64501} & \num{-13.63207} & \num{54.09503} 
                    &   \num{-23.56024} & \num{-13.60452} & \num{53.47486} 
                    &   \num{-23.55939} & \num{-13.60425} & \num{53.46866} \\
        \bottomrule
    \end{tabular}
\end{table}

As for the case of pure bending, we obtain accurate results for the present problem, which is characterized by a more complex state of strain. 
Even for coarse meshes and low element orders, tip-displacements $\uv(L)$ are within less than one percent of the reference solution $\uv_{\mathrm{ref}}$, which is computed with $512$ fourth-order ($k=4$) elements.
The converged results for $\rho=100$ coincide exactly with those reported by Meier et al.~\cite{Meier2019}\footnote{
    Note that, though stated otherwise, tip positions are given in~\cite{Meier2019} rather than displacements.
}.
To quantify the rates of convergence, we compute the $L^2$-error in the displacements, which is normalized by the radius $R$ of the arc,
\begin{equation}
    \err = \frac{1}{R} \sqrt{\frac 1 L \int_0^L \left\Vert \uv - \uv_{\mathrm{ref}}\right\Vert^2 ds} .
\end{equation} 
\begin{figure}[htb!]
    \centering
    \includegraphics[width=1.\textwidth]{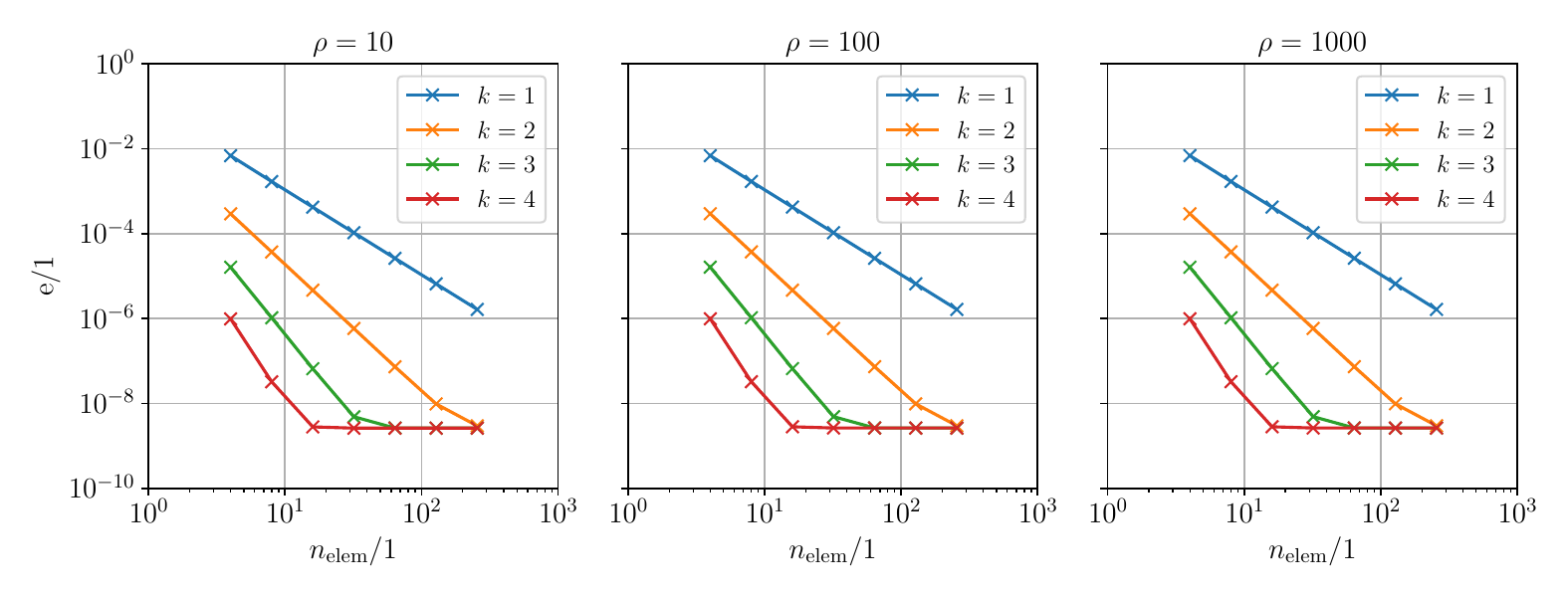}
    \caption{Arc-segment under out-of-plane force: rates of convergence for beam elements of orders $k=1,\ldots,4$ and different values of slenderness; $\rho = 10$ (left), $\rho = 100$ (center), and $\rho = 1000$ (right).}
    \label{fig:3-arc_convergence}
\end{figure}

The convergence plots illustrated in Fig.~\ref{fig:3-arc_convergence} reveal optimal rates of convergence for the considered element orders.
We note that the slenderness does not affect the convergence in any way. 
As a matter of fact, the errors are visually indistinguishable for the three values of the slenderness considered.
The minimum attainable normalized error is of the order $10^{-9}$ for the present problem irrespective of the beam's slenderness.

\subsection{Helix experiment}
In this benchmark example, which is taken from \cite{harsch2021}, an initially straight beam is deformed to take the shape of a perfect helix. While position and rotation are prescribed accordingly at one end, an external moment load is applied to the other end of the beam. The helix has $n=2$ coils of radius $R_0=10$ over a height $h_0=50$. 
Under these assumptions, the total length of the beam amounts to $L = \sqrt{1+c^2}\left(2 \pi n R_0 \right)$ for $c = h_0/\left(2\pi n R_0\right)$, and the beam's center line position is given by
\begin{align}
    \rv (s) = R_0 \sin(\alpha(s)) \ev_1 - R_0 \cos(\alpha(s)) \ev_2 + c R_0 \alpha(s) \ev_3, \qquad \alpha(s) = 2\pi n s / L.
    \label{eq:helix}
\end{align}
Based on the slenderness ratio $\rho = 10$, for a cylindrical beam we obtain a cross-section radius $r = L/(2\rho)$, cross-sectional area $A = r^2 \pi$ and moments of area $I_2 = I_3 = r^4\pi/4$, $I_1 = r^4\pi/2$. A Young's modulus $E=1$ and shear modulus $G=0.5$ lead to axial stiffness $\EA$, shear stiffness $\GAy=\GAz=GA$, as well as equivalent bending and torsional stiffnesses $\EIy = \EIz=\GIt$.
For details on the derivation of the respective choices for boundary conditions we refer to the original publication by Harsch et al.~\cite{harsch2021}. Under the condition of equivalent bending and torsional stiffnesses, the following boundary conditions lead to deformation state \eqref{eq:helix},
\begin{align}
    \rv_0(0) &= -R_0 \ev_2, & \psiv(0) &= \arctan{(c)}\, \ev_2,\\
    \nv(L) &= \boldsymbol{0} ,& \mv(L) &= \frac{1}{R_0 \left(1+c^2\right)} \left( \GIt c\, \gv_1 + \EIz\, \gv_3\right).
\end{align}

\begin{figure}
    \centering
    \includegraphics[width=0.9\textwidth]{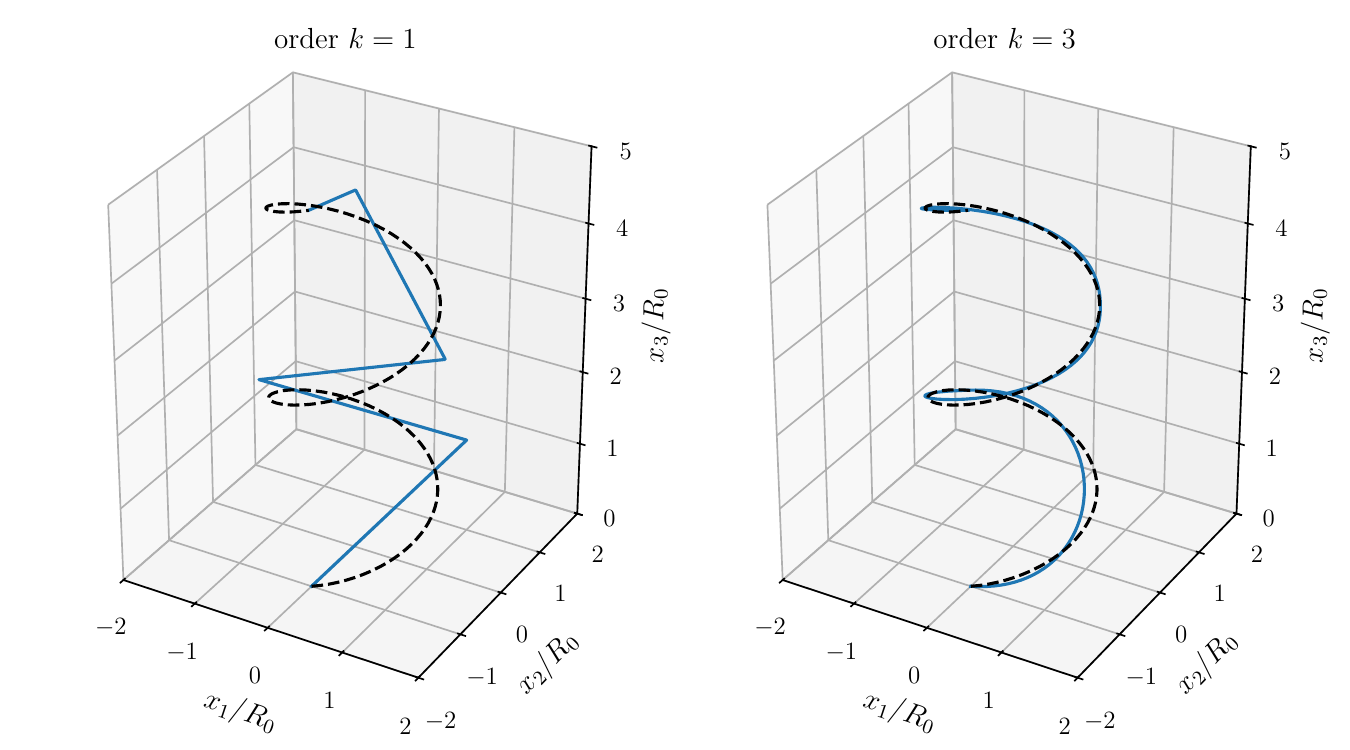}\\
    \includegraphics[width=0.9\textwidth]{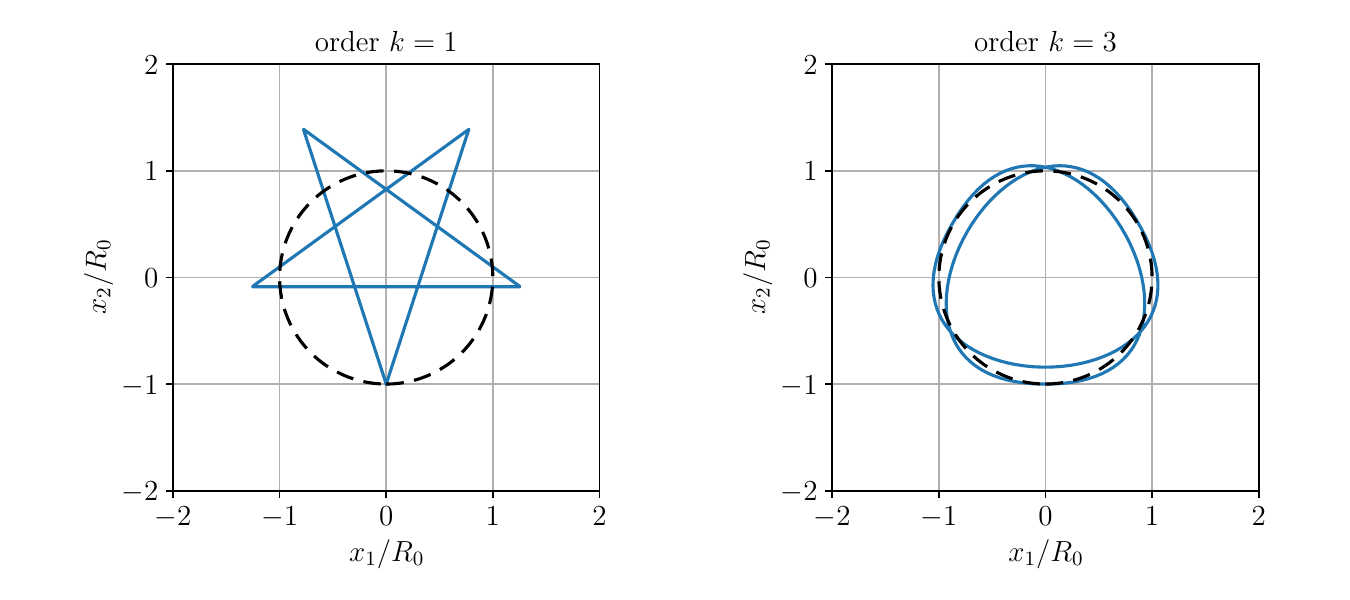} 
    \caption{Helix experiment: computed configurations for very coarse discretizations of $5$ elements of order $k=1$ (left) and $3$ elements of order $k=3$ (right), respectively. The position of the tip point is obtained numerically exact in both computations.}
    \label{fig:helixshapes}
\end{figure}

The position of the tip node is obtained exactly (i.e., up to the floating point accuracy) for as few as $5$ elements of lowest order. Unlike the $SE(3)$ element presented in \cite{harsch2023}, the current element is not restricted to relative rotations bounded by $\pi$ throughout the element. When using at least second-order elements, a total of $3$ elements is enough to reproduce numerically exact nodal values for the beam position and rotation. However, relative rotations of more than $2\pi$ cannot be represented within a single element. In Fig.~\ref{fig:helixshapes}, the computed center lines for very coarse discretizations of five first-order and three third-order elements are shown. 

\subsection{Rod bent to helical shape}
Ibrahimbegovic~\cite{ibrahimbegovic1997a} proposed the following example to explore the capabilities of beam elements regarding the representation of large, inhomogeneous deformations.
A concentrated moment and force are applied in simultaneous load steps to the tip of an initially straight cantilever beam, such that a helical shape is obtained. 
A beam of length $L=10$, axial stiffness $\EA = 10^4$, shear stiffnesses $\GAy=\GAz = \num{1e4}$, bending stiffnesses $\EIy = \EIz = \num{1e2}$ and torsional stiffness $\GIt = \num{1e2}$ is considered; an incremental application of a total tip moment $\mv = 20\pi \EIy/L \ev_3$ and out-of-plane tip load $\nv = 50\ev_3$ induces a helical shape mounting up to slightly over ten full circles of the beam. To this end, a static load factor $\lambda \in [0,1]$ is introduced which appears in the virtual work of external forces,
\begin{align}
    \delta W_{ext} = \left.\left(\lambda \nv \cdot \delta \rv + \lambda \mv \cdot \delta \thetav\right)\right|_{s=L}.
\end{align}

\begin{figure}[htb!]
    \centering
    \includegraphics[width=0.5\textwidth]{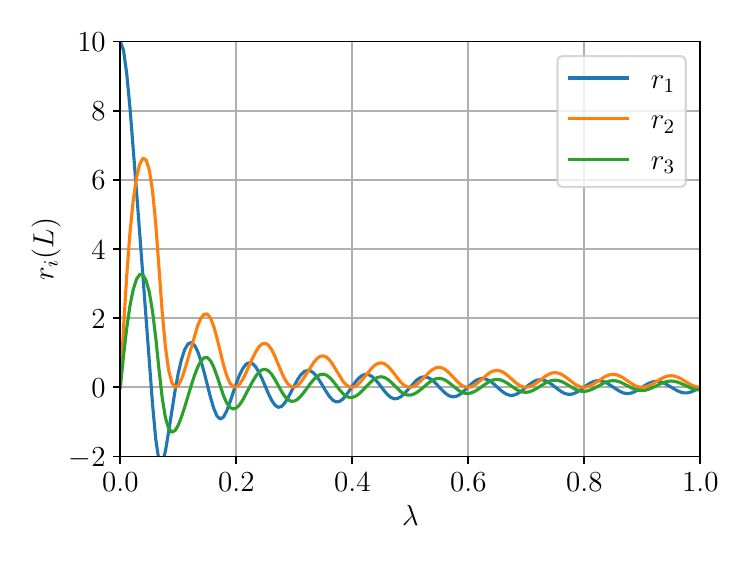}
    \caption{Rod bent to a helical shape: evolution of tip position over static load factor $\lambda$ for a discretization of $30$ elements of polynomial orders $k=3$.}
    \label{fig:helical_tipdisp}
\end{figure}

The out-of-plane component of the tip position, and thereby the observed height of the helical shape, does not grow continuously with $\lambda$, but shows oscillatory behavior. In Fig.~\ref{fig:helical_tipdisp}, the evolution of all three components of the tip displacement is provided. The results are in direct agreement with the original findings by Ibrahimbegovic~\cite{ibrahimbegovic1997a}, as well as Harsch et al.~\cite{harsch2023}. 

\begin{figure}[htb!]
    \centering
    \includegraphics[width=0.95\textwidth]{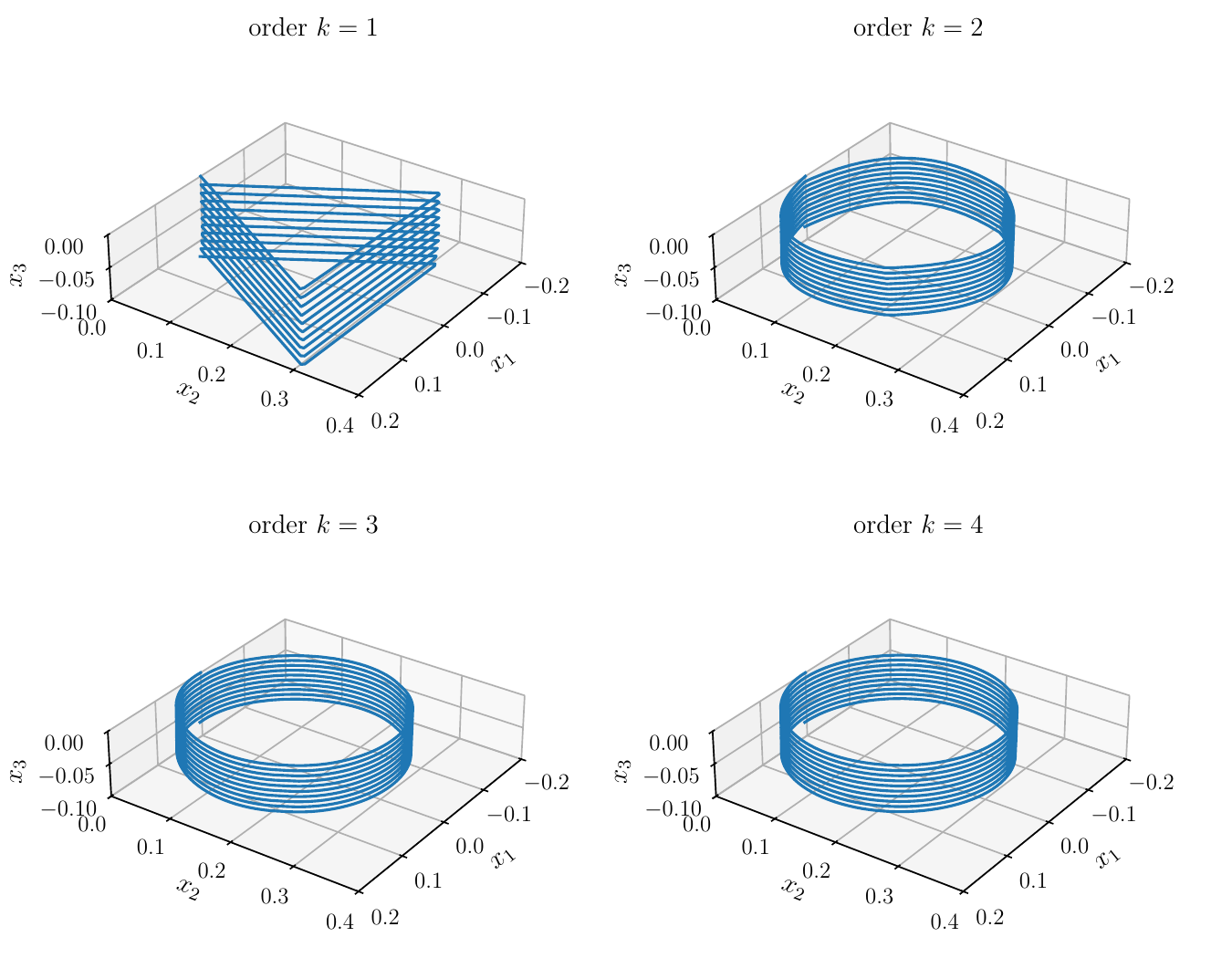}
    \caption{Rod bent to a helical shape: final configuration for a discretization of $30$ elements of polynomial orders $k=1,\ldots,4$.}
    \label{fig:helicalshapes}
\end{figure}

To assess the accuracy of the method, the tip displacement computed on a discretization of $30$~elements of polynomial order one to four is compared to the value obtained on a finer mesh consisting of $100$ fourth-order elements. Figure~\ref{fig:helicalshapes} shows the computed beam center line in final configuration for these various polynomial orders; Table~\ref{tab:helical_tipdisp} provides numerical values for the tip displacement. Even for first-order elements,  the computation remains stable up to $\lambda=1$, although one circle is resolved by only three elements. Given this coarse discretization, the accuracy of the solution is still satisfactory. As expected, higher-order elements yield more accurate results.

\begin{table}[htb!]
    \caption{Rod bent to a helical shape: components of tip position $r_1$, $r_2$, $r_3$ observed in final configuration for a discretization of $30$ elements of polynomial orders one to four, and component-wise absolute error $e_1$, $e_2$, $e_3$ as compared to solution for $100$ elements of order five.}
    \centering
    \begin{tabular}{c|ccc|ccc}
        \toprule
        $k$ & $r_1(L)$ & $r_2(L)$ & $r_3(L)$ & $e_1$ & $e_2$ & $e_3$ \\
        \midrule
        $1$   & \num{6.698e-03} &   \num{1.166e-04} & \num{-1.249e-1} &
        \num{1.9e-03} &  \num{4.4e-05} &  \num{4.8e-02}\\ 
        $2$   & \num{4.753e-03} &  \num{7.135e-05} & \num{-7.789e-02} &
        \num{5.0e-05} & \num{1.2e-06} & \num{1.5e-03}\\
        $3$   & \num{4.803e-03} & \num{7.248e-05} & \num{-7.647e-02} &
        \num{5.6e-07} & \num{2.0e-08} & \num{2.4e-05} \\
        $4$   &  \num{4.803e-03} & \num{7.250e-05} & \num{-7.644e-02} &
        \num{9.8e-09} & \num{2.8e-10} & \num{8.7e-08} \\
        \bottomrule
    \end{tabular}
    \label{tab:helical_tipdisp}
\end{table}

\begin{figure}
    \centering
    \includegraphics[width=0.95\textwidth]{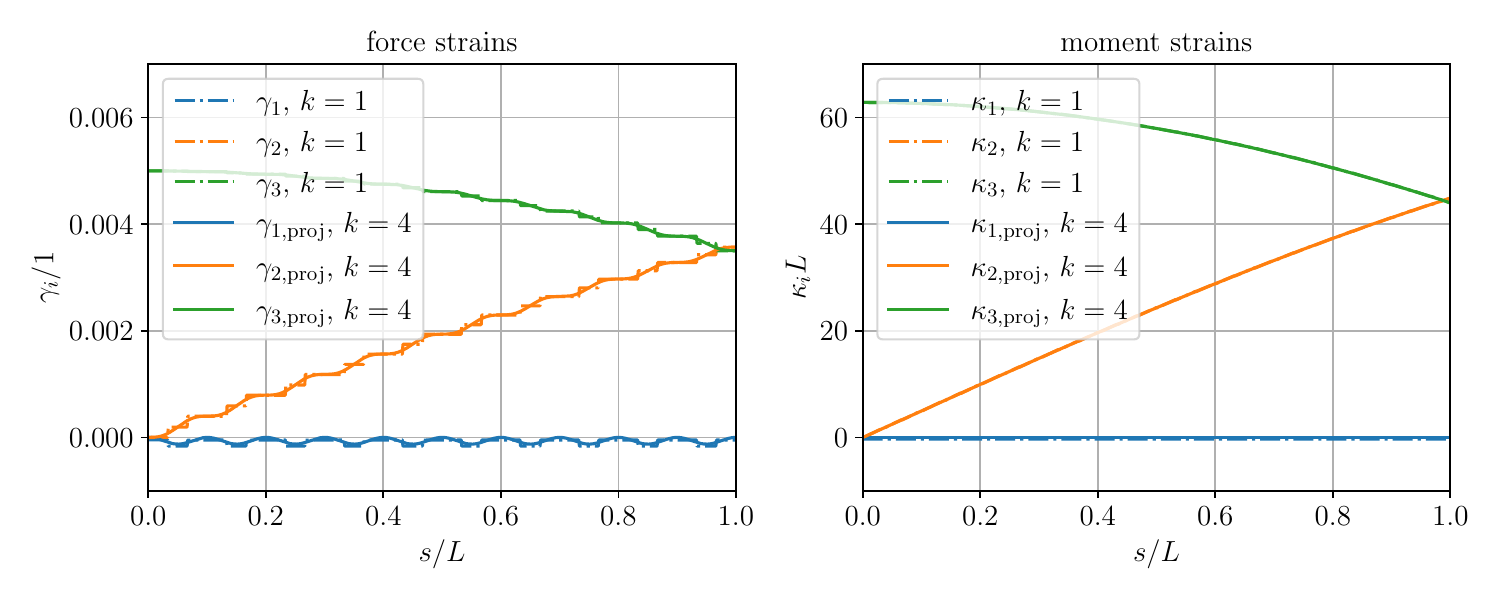}
    \caption{Rod bent to a helical shape: force and moment strains over axial coordinate $s$ for a discretization of $30$ elements and polynomial orders $k=1$ and $k=4$.}
    \label{fig:helical_strains}
\end{figure}

Harsch et al.~\cite{harsch2023} also provide an evaluation of strain measures $\gammav$ and $\kappav$ along the axis at final configuration. Figure~\ref{fig:helical_strains} shows strain and moment strain along the beam's axis for the discretizations of order $k=1$ and $k=4$. Note that, for higher-order elements, the force strain $\gammav$ is projected element-wise to polynomials of degree $k-1$. This procedure is in line with the application of reduced integration described in Section~\ref{sec:locking}, and yields generalized strains of optimal order given displacement discretization of order $k$. To evaluate moment strains $\kappav$, the inverse material relations are employed, evaluating $\kappav = (\Ckappa)^{-1} \cdot \Mv$. For this choice of representation, moment strains are piecewise of order $k$.

\subsection{Discontinuous slope}

\begin{figure}[htb!]
    \centering
    \includegraphics[width=0.45\textwidth]{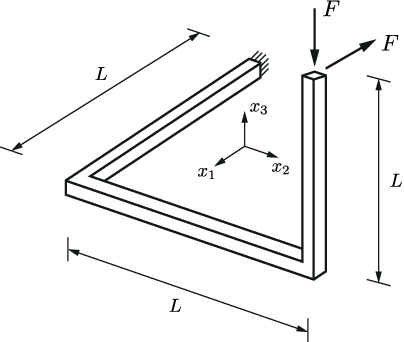}
    \hfill
    \includegraphics[width=0.5\textwidth]{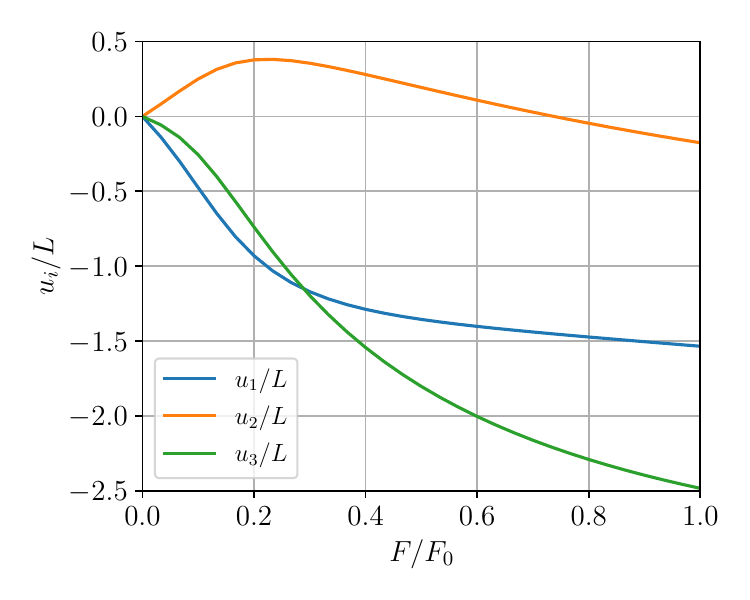}
    \caption{Structure with slope discontinuity: Geometry and external loads (left, reproduced from~\cite{romero2008}); components of the tip-displacement as a function of the load intensity (right).}
    \label{fig:5_slope_discont_geo}
\end{figure}

The proposed element allows us to analyze structures with kinks and branches without much ado.
As a first case, we consider the example of three straight rods connected at right angles, which was originally proposed by Romero~\cite{romero2008} and later studied, e.g., in~\cite{eugster2014,Wasmer2024}, see Fig.~\ref{fig:5_slope_discont_geo} (left) for an illustration of the geometry.
Each member of the structure has unit length $L = 1$.
The beam is clamped at $s=0$; a concentrated force $\nv(3L) = - F \ex - F \ez$ is applied at the free end.
In line with~\cite{Wasmer2024}, the effective stiffnesses are chosen as $\EA = 2 \GAy = 2 \GAz = \num{1e4}$, $\GIt = \EIy = \EIz = 100 / 12$.

\begin{figure}[htb!]
    \centering
    \includegraphics[width=\textwidth]{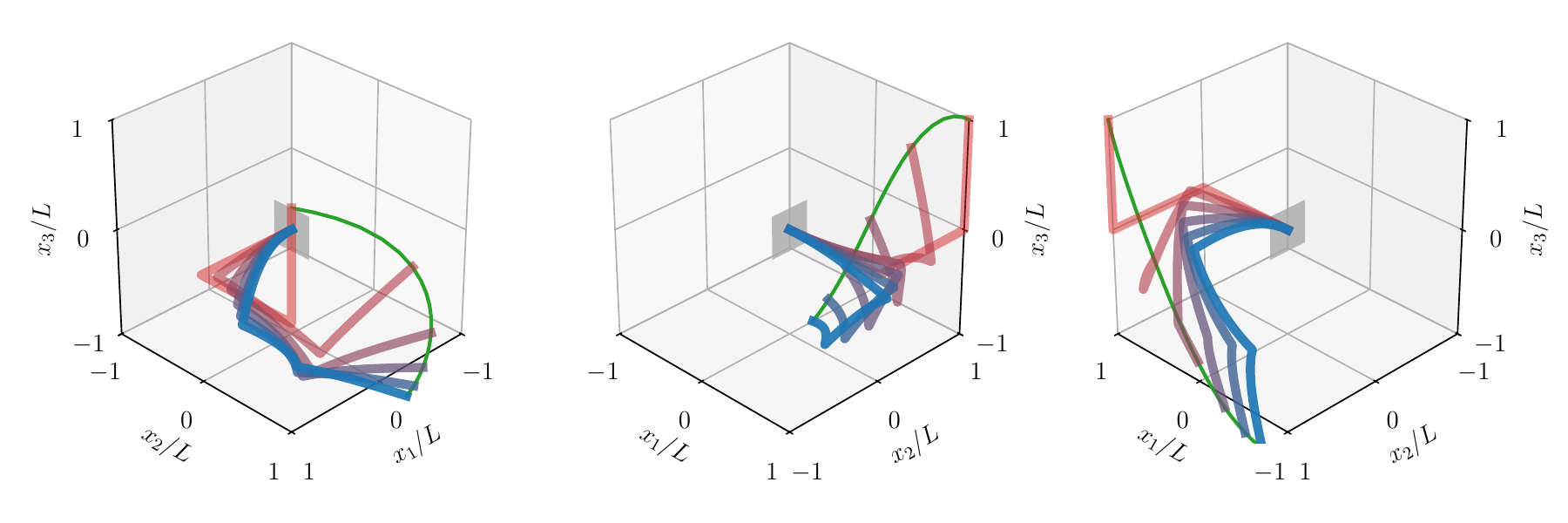}
    \caption{Structure with slope discontinuity: series of configurations for loads $F = 0,\ldots,10$ in increments of 2. The green line represents the trajectory $\rv(3L)$ of the loaded tip.}
    \label{fig:5_slope_discont_configs}
\end{figure}

Figure~\ref{fig:5_slope_discont_configs} shows a series of configurations of the structure from different perspectives, where the load intensity is raised up to a maximum value of $F_0 = \num{10}$ in increments of \num{2}. 
Colors are blended from red (undeformed configuration) to blue (deformed configuration) for the sake of clarity; the green line indicates the trajectory of the free end ($s=3L$) in the course of loading.
The components of tip-displacement are illustrated as function of the load intensity in Fig.~\ref{fig:5_slope_discont_geo} (right).

The values of the tip displacement in the fully loaded configuration are listed in Tab.~\ref{tab:5_slope_discont_tip_displ} for element orders $k=1,\ldots,4$.
The top half provides the results for the coarsest possible discretization, where each member is represented by a single element only. 
The results in the lower half of the table are obtained with eight elements per segment, i.e., $\nel = \num{24}$ elements in total.
As with the previous examples, the effectiveness of the proposed element is also confirmed for the present problem. 
A single fourth-order element per segment is sufficient to obtain converged results in spite of the general state of (large) deformation.
With eight elements per segment, results are converged for orders $k \geq 2$, but even the lowest-order element ($k=1$) provides an excellent approximation showing sub-permille error in the tip-displacement.

\begin{table}[htb!]
    \sisetup{
    round-mode = figures,
    round-precision = 6,
    group-digits = integer
    }
    \centering
    \caption{Structure with slope discontinuity: components of the tip-displacement $\uv(L)$ for  element orders $k=1,\ldots,4$ obtained with one element per segment (top) and eight elements per segment (bottom).}
    \label{tab:5_slope_discont_tip_displ}
    \begin{tabular}{c c c c c}
        \toprule
        $k$ &   $\nel$ &   $u_1/L$    &   $u_2/L$    &   $u_3/L$ \\
        \midrule
        1   &   3   &   \num{-1.519637} & \num{-0.1870366} & \num{-2.607832} \\
        2   &   3   &   \num{-1.528034} & \num{-0.1713543} & \num{-2.470063} \\
        3   &   3   &   \num{-1.535102} & \num{-0.175832} & \num{-2.484094} \\
        4   &   3   &   \num{-1.535068} & \num{-0.1758733} & \num{-2.484223} \\
        \midrule
        1   &   24   &   \num{-1.534489} & \num{-0.1759337} & \num{-2.485712} \\
        2   &   24   &   \num{-1.53507} & \num{-0.1758742} & \num{-2.484215} \\
        3   &   24   &   \num{-1.535072} & \num{-0.1758755} & \num{-2.484219} \\
        4   &   24   &   \num{-1.535072} & \num{-0.1758755} & \num{-2.484219} \\
        \bottomrule
    \end{tabular}
\end{table}

\subsection{Fork structure}
\begin{figure}[htb!]
    \centering
    \includegraphics[width=0.45\textwidth]{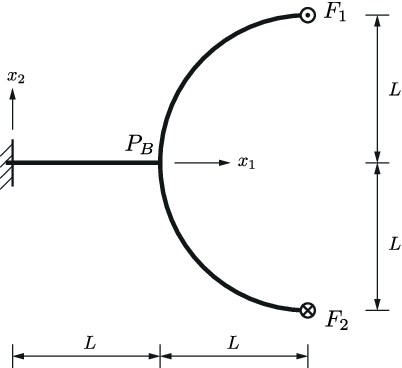}
    \caption{Fork structure subjected to a force couple: undeformed geometry and external loads.}
    \label{fig:7_fork}
\end{figure}

As final example, we consider a fork-like structure as illustrated in Fig.~\ref{fig:7_fork}.
It is composed from a straight section of length $L = \num{1}$, which branches into semi-circle with a radius of $L$.
The fork, which lies in the $x_1x_2$-plane in the undeformed configuration, is subjected to a pair of opposing out-of-plane forces applied at the tips of the fork's tines, i.e., the force $\Fv_1 = F \ez$ is acting at the point $P_1 = (2L,L,0)$, whereas $\Fv_2 = -F \ez$ is applied at $P_2 = (2L,-L,0)$.
The branching point is denoted by $P_B = (L,0,0)$ subsequently.
The cross-sectional stiffnesses are set to $\EA = \GAy = \GAz = \num{1e4}$ and $\GIt = \EIy = \EIz = \num{1e2}$, respectively.
The purpose of the present problem is twofold. By means of the example, we can show that the proposed formulation can be immediately applied to structures with branching centerlines, since the concept of discrete curvature allows us to deal with branching points in a very straightforward way.
Secondly, the presence of a pair of forces enables us to numerically verify the path-independence of the beam element, which we expect from the fact that the formulation forgoes deformation increments but involves total displacement and rotation fields only.
To this end, the forces are applied consecutively. 
Assuming branches (no pun intended) in the deformation path to be absent, the same equilibrium configuration needs to be recovered irrespective of the order of loading.
As a matter of fact, the nature of the present problem allows us to demonstrate the complete symmetry of the solution upon reversal of the order.
%

\begin{figure}[htb!]
    \centering
    \includegraphics[width=0.32\textwidth,trim=0.8cm 1mm 2cm 0,clip=true]{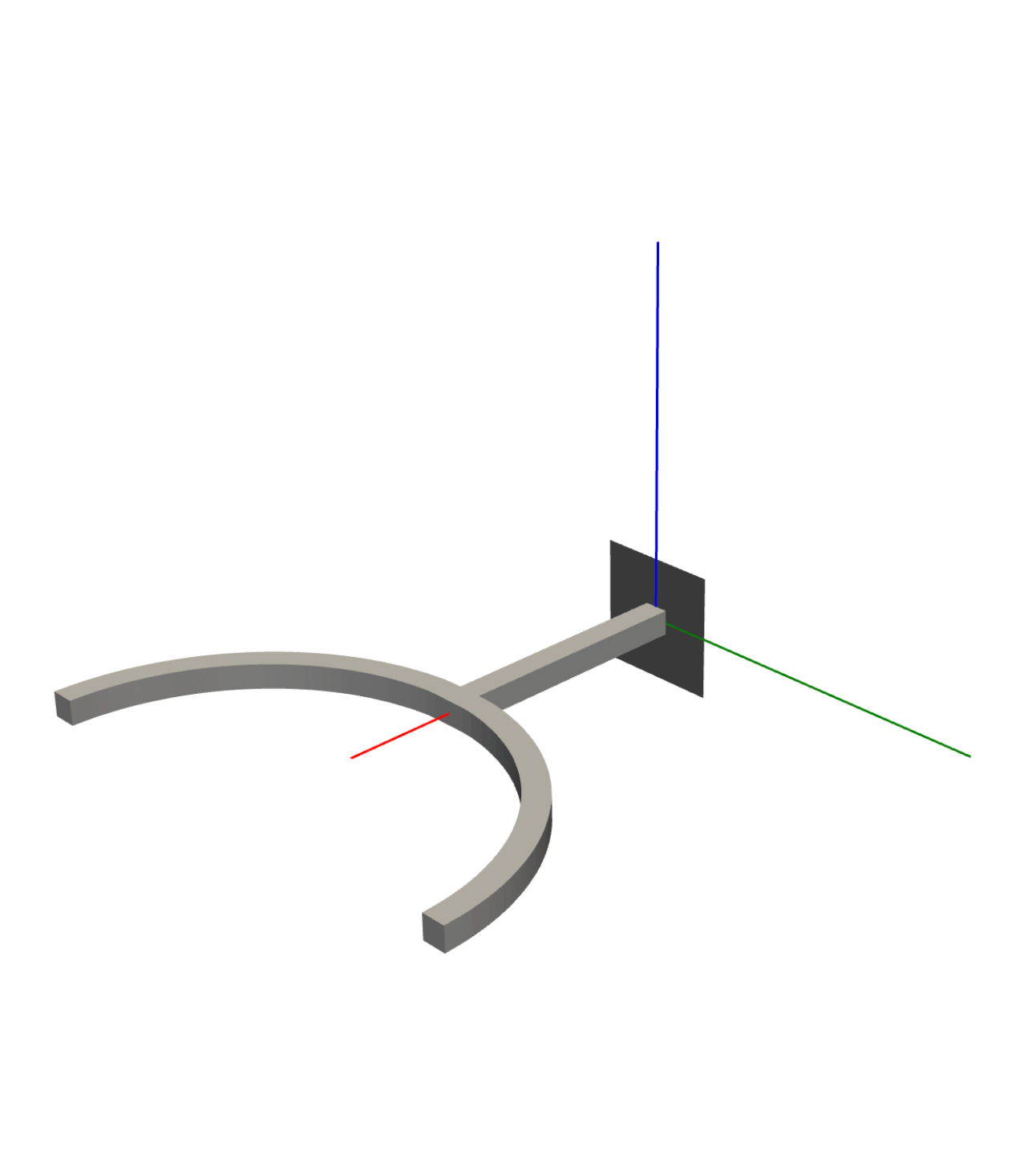}
    \includegraphics[width=0.32\textwidth,trim=0.8cm 1mm 2cm 0,clip=true]{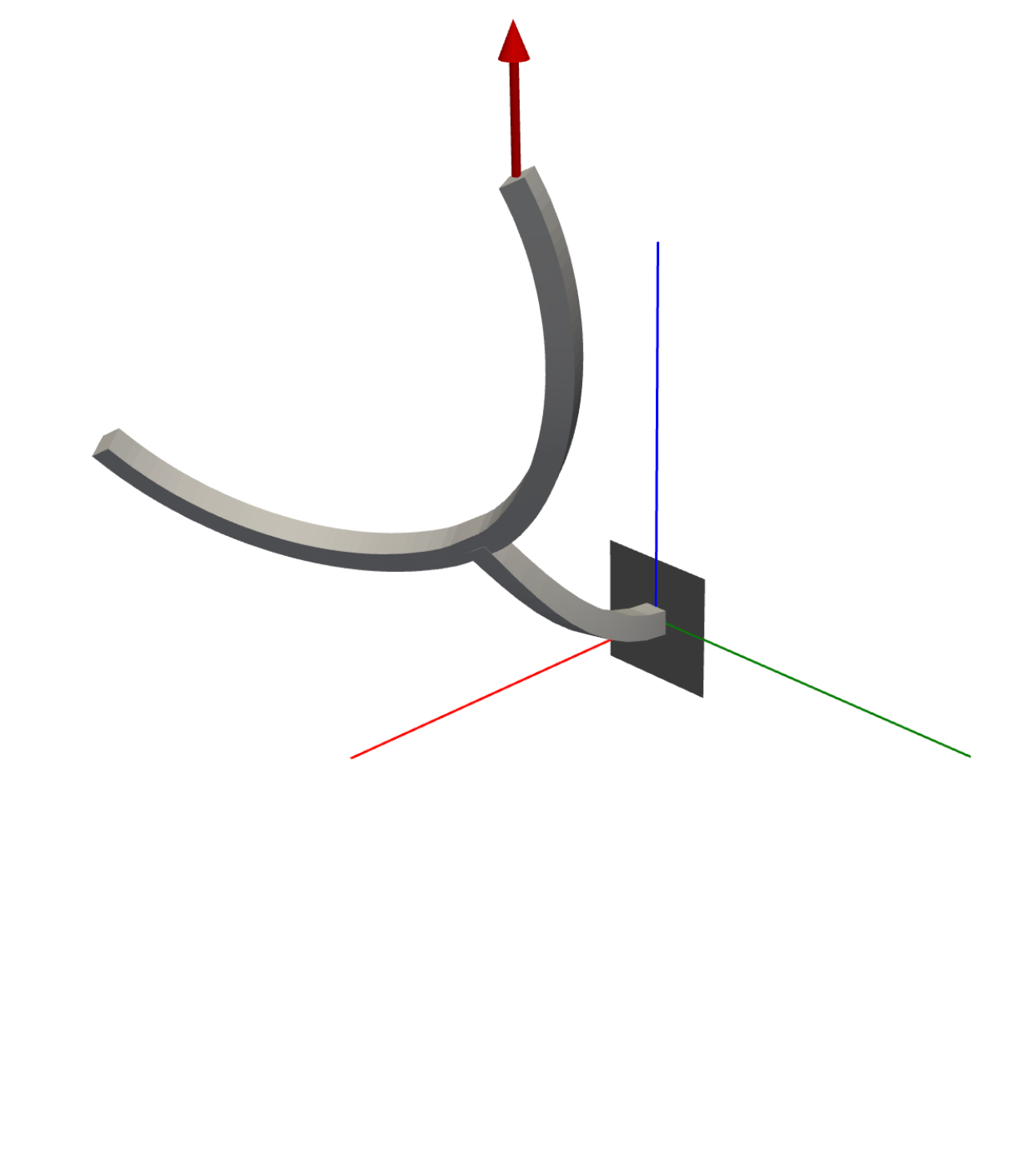}
    \includegraphics[width=0.32\textwidth,trim=0.8cm 1mm 2cm 0,clip=true]{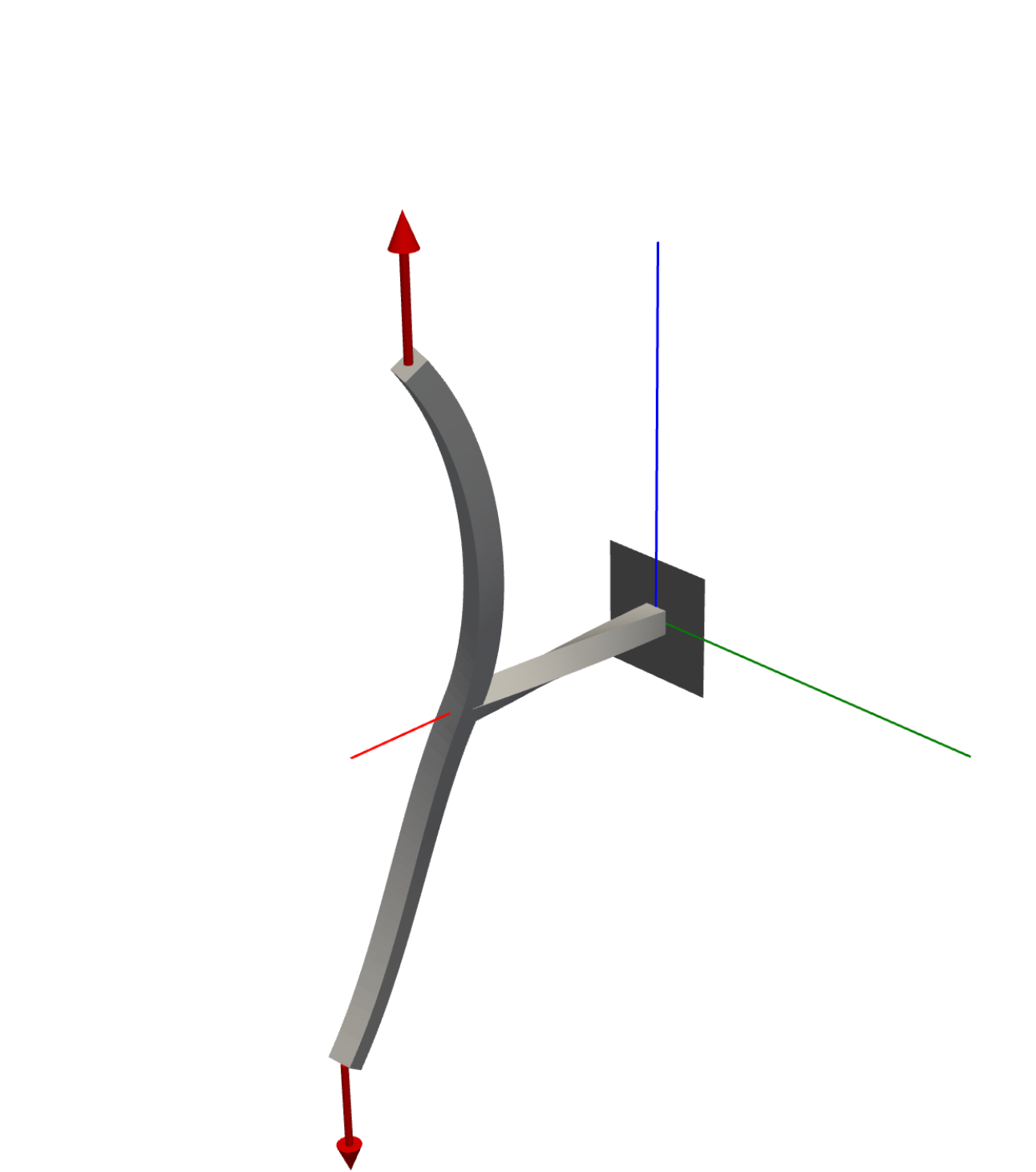}
    \caption{Fork structure subjected to a force couple: undeformed configuration (0; left), equilibrium configuration under a single force (1; center), and final configuration under a force couple (2; right). Red arrows indicate the forces acting in the respective configurations.}
    \label{fig:7-fork_configs}
\end{figure}

Starting from the unloaded configuration $(0)$, we first increase the force applied at $P_1$ up to a nominal value of $F = \num{200}$ over ten load increments, keeping the other tine ($P_2$) of the fork unloaded.
From the obtained equilibrium configuration $(1)$, the force at $P_2$ is imposed over ten more load steps, after which we arrive at the final configuration $(2)$, in which the force couple has fully developed. 
Figure~\ref{fig:7-fork_configs} illustrates the undeformed configuration (left), the equilibrium configuration under a single force $\Fv_1$ (center), and the final configuration, in which both $\Fv_1$ and $\Fv_2$ are applied (right).
When subjected to a single force, the fork obliquely bends out of the $x_1x_2$-plane due to the non-symmetric load.
After imposing the second opposing force, the shaft of the fork returns to the $x_1$-axis and is in a pure state of twist, while the two tines are significantly bent apart.

\begin{figure}[htb!]
    \centering
    \includegraphics[width=\textwidth]{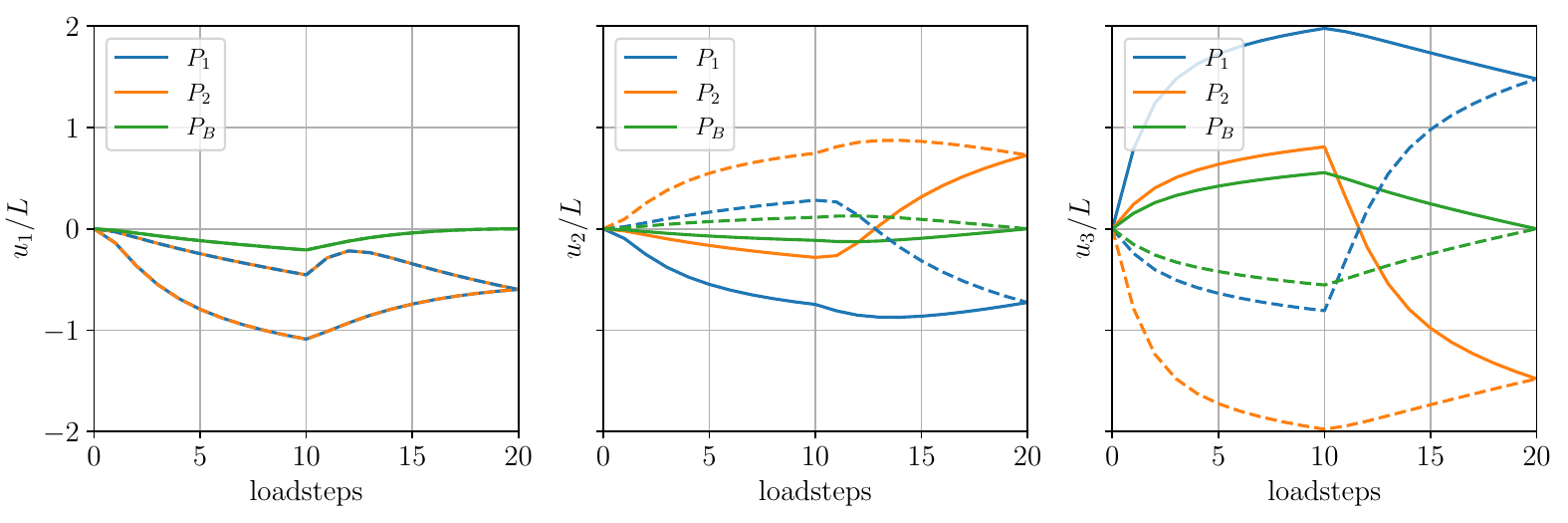}
    \caption{Fork structure subjected to a force couple: components of the displacements of the two tines $P_1$, $P_2$ and the branching point $P_B$. Dashed lines correspond to a reversal of the load application.}
    \label{fig:7-fork_displ}
\end{figure}

In Fig.~\ref{fig:7-fork_displ}, the components of the displacement of the two tines ($P_1$, $P_2$), and the displacement of the branching point ($P_B$) are illustrated over the load steps. 
The dashed lines correspond to a reversed loading, i.e., the force at $P_2$ is applied over the first ten load steps, before the force $P_1$ is imposed.
The obtained load-displacement curves reflect the symmetry of the fork's geometry and loading, lacking any indication of path-dependence of the proposed formulation. 

Numerical values of the displacement vector of the first tine's tip $P_1$ in configurations $(1)$ and $(2)$ are listed in Tab.~\ref{tab:7-fork_displ}. 
As with the previous examples, we study element orders $k=1,\ldots,4$ and two levels of mesh refinement.
The coarse mesh (upper half) uses three elements per segment, whereas ten elements are employed in the fine mesh (lower half), which translates into a total of $\nel = \num{9}$ and $\nel = \num{30}$, respectively.
Once again, we can observe excellent results, with the tip-displacement being already converged for the coarse mesh when using fourth-order elements ($k=4$).
For the fine mesh, even linear elements ($k=1$) show relative errors in the displacements on a per-mille level.

\begin{table}[htb!]
    \sisetup{
    round-mode = figures,
    round-precision = 6,
    group-digits = integer
    }
    \centering
    \caption{Fork structure: components of the tip-displacement $\uv$ at tine $P_1$ in configurations (1) and (2) for element orders $k=1,\ldots,4$ obtained with three elements per segment (top) and ten elements per segment (bottom).}
    \label{tab:7-fork_displ}
    \begin{tabular}{c c | c c c | c c c }
        \toprule
        $k$ &   $\nel$  &   $u_1^{(1)}/L$    &   $u_2^{(1)}/L$    &   $u_3^{(1)}/L$
                        &   $u_1^{(2)}/L$    &   $u_2^{(2)}/L$    &   $u_3^{(2)}/L$ \\
        \midrule
        1   &   9   &   \num{-1.084348} & \num{-0.759430} & \num{1.981360} 
                    &   \num{-0.599783} & \num{-0.731959} & \num{1.465554} \\
        2   &   9   &   \num{-1.088565} & \num{-0.746239} & \num{1.978099} 
                    &   \num{-0.598926} & \num{-0.727164} & \num{1.479676} \\
        3   &   9   &   \num{-1.088609} & \num{-0.746318} & \num{1.978403} 
                    &   \num{-0.598811} & \num{-0.727135} & \num{1.479965} \\
        4   &   9   &   \num{-1.088614} & \num{-0.746307} & \num{1.978303} 
                    &   \num{-0.598829} & \num{-0.727134} & \num{1.479862} \\
        \midrule
        1   &   30  &   \num{-1.088214} & \num{-0.747378} & \num{1.978359} 
                    &   \num{-0.598936} & \num{-0.727541} & \num{1.478570} \\
        2   &   30  &   \num{-1.088613} & \num{-0.746306} & \num{1.978299} 
                    &   \num{-0.598830} & \num{-0.727134} & \num{1.479858} \\
        3   &   30  &   \num{-1.088614} & \num{-0.746307} & \num{1.978302} 
                    &   \num{-0.598829} & \num{-0.727134} & \num{1.479861} \\
        4   &   30  &   \num{-1.088614} & \num{-0.746307} & \num{1.978301} 
                    &   \num{-0.598829} & \num{-0.727134} & \num{1.479860}\\
        \bottomrule
    \end{tabular}
\end{table}

\section{Conclusion}
\label{sec:conlusion}

Our work shows that, in the field of nonlinear beam theories---although it is by no means a new topic---there are still things to discover. 
Specifically, we have succeeded in deriving a new mixed method for geometrically exact beams, and, on this basis, we have developed a beam element that shows superior numerical properties.
For this purpose, the moment vector has been introduced as an independent field, complementing positions and rotations.
A discontinuous interpolation of the rotations has proven central to our method. 
For the virtual work to be well-defined at element interfaces (vertices), the notion of a discrete curvature has been introduced.
We have shown that the moment vector is work-conjugate to the change in relative rotation across the interface of adjacent elements.
The discontinuous interpolation of rotations allows us to establish objectivity of the discretized strain measures (almost) en passant.
We have introduced a multiplicative split of the element-local rotation vector into a constant lowest-order part and higher-order relative rotations, from which curvatures are computed.
Hybridizing the element---i.e., restricting the continuity of moments to the element interior and adding nodal rotations as additional unknowns---allows conventional boundary conditions to be imposed in a straightforward way.

The effectiveness of the proposed formulation has been shown with several classical benchmark problems.
Objectivity and path-independence have been demonstrated; when combined with reduced integration in the strain energy contribution associated with generalized force strains, our beam element shows optimal rates of convergence irrespective of the polynomial order employed and the beam's slenderness.
Contained in our formulation is a lowest-order element that entirely dispenses with the interpolation of rotations, which is particularly appealing in view of the expensive arithmetic involved in finite rotations.

Summarizing, the numerical examples reveal the great potential of the proposed mixed formulation and its realization as beam element.
Still, there are open questions that remain to be answered. 
For the moment, we have only considered the static case. Dynamic problems along with appropriate time-integration schemes will be subject to future work.
It is certainly appealing to investigate whether and how concepts of our approach can be transferred to the motion formalism, which aims at a discretization of $SE(3)$.
Beams just never seem to get old.

\bibliography{refs}

\end{document}